\documentstyle[12pt]{article}
\newcommand{\1}{{\rm 1\!\!1}}

\newcounter{mysection}
\def\myownsection{\refstepcounter{mysection} \setcounter{equation}{0}}

\begin{document}
$\;$\\[20pt]
\begin{center}
{\bf FILTERED RANDOM VARIABLES,\\
BIALGEBRAS AND CONVOLUTIONS}\\[60pt]
{\sc Romuald Lenczewski}\\[40pt]
Institute of Mathematics\\ 
Wroc{\l}aw University of Technology\\
Wybrzeze Wyspianskiego 27\\
50-370 Wroc{\l}aw, Poland\\
e-mail lenczew@im.pwr.wroc.pl\\[40pt]
\end{center}
\begin{abstract}
We introduce the {\it filtered *-bialgebra} which is
a multivariate generalization of the unital *-bialgebra
${\bf C}\langle X, X', P\rangle $ of polynomials in noncommuting
variables
$X=X^{*}$, $X'^{*}=X'$ and a projection $P=P^{*}=P^{2}$, 
endowed with the coproduct
$\Delta (X)=X\otimes 1 + 1\otimes X$, 
$\Delta (X')=X'\otimes P + P \otimes X'$, with $P$ being group-like. 
We study the associated convolutions, random walks
and {\it filtered random variables}.
The GNS representations of the limit states
lead to {\it filtered fundamental operators} 
which are the CCR fundamental operators
on the multiple symmetric Fock space $\Gamma({\cal H})$ over
${\cal H}=L^{2}({\bf R}^{+},{\cal G})$, where ${\cal G}$ is a
separable Hilbert space, multiplied by appropriate projections.
The importance of filtered random variables and fundamental operators
stems from the fact that by addition and strong limits one obtains from them
the main types of 
noncommutative random variables and fundamental operators, respectively,
regardless of the type of noncommutative independence.
\\[5pt]
Mathematics Subject Classification (1991): 81R50, 60J15, 46L50\\[10pt]
\end{abstract}
\myownsection
\begin{center}
{\sc 1. Introduction}
\end{center}
In this work we introduce and study basic noncommutative
random variables, from which the main types of noncommutative
random variables can be constructed regardless of the notion of
independence.

The basic idea of introducing filtered random variables is 
pretty straightforward and has its origin in the definition of the 
convolution of measures and the associated states.
Let ${\bf C}[X]$ be the unital *-algebra of polynomials in $X^{*}=X$, with 
the coproduct 
\begin{equation}
\label{1.1}
\Delta (X)=X\otimes 1 + 1 \otimes X.
\end{equation}
If $\phi, \psi$ are states on ${\bf C}[X]$, then 
$$
\phi \star_{c} \psi =\phi \otimes \psi \circ \Delta
$$
gives the convolution of states corresponding to the classical convolution
of measures. 

In order to define
a quantum deformation of this simple model, we replace the unit
in the coproduct (1.1) by a projection $P$ to get
\begin{equation}
\label{1.2}
\Delta (X')=X'\otimes P + P \otimes X', \;\;\;
\Delta (P)=P\otimes P
\end{equation}
Then, for given state $\phi$ on 
${\bf C}[X']$, we define its noncommutative extension
$\widetilde{\phi}$ to ${\bf C}\langle X', P \rangle$, which
is the free product ${\bf C}[X']*{\bf C}[P]$ with identified units, by
$$
\widetilde{\phi}(P^{\alpha}Y^{n_1}PY^{n_2}P\ldots Y^{n_k}P^{\beta})
=\phi(Y^{n_1})\phi(Y^{n_2})\ldots \phi(Y^{n_k})
$$
where $\alpha , \beta \in \{0,1\}$ and $n_{1}, \ldots , n_{k}\in {\bf N}$,
called the Boolean extension [11].
The convolution 
$$
\widetilde{\phi} \star_{B}  \widetilde{\psi}
=\widetilde{\phi}\otimes \widetilde{\psi}\circ \Delta
$$
where $\Delta $ is given by (1.2), 
gives a quantum analog of the classical convolution of states, 
called the Boolean convolution. Note that by introducing $P$ we can
deal with the usual tensor coproducts instead of the special one as in
the approach of Sch\"{u}rmann [19]. 
The same holds for the $m$-free and free products of states.

This new convolution is very important since its generalization
to the multivariate case, when restricted to 
suitable *-bialgebras, gives also $m$-free and 
free convolutions (see [4] and [11]).
In the multivariate case we study the unital *-algebra
$\widehat{\cal B}$
over ${\bf C}$ generated by $X_{k}(\sigma), P(\sigma)$
$k\in {\bf N}, \sigma\in {\cal P}({\bf N})$, where
${\cal P}({\bf N})$ is the power set of ${\bf N}$,
with the involution given by
$X_{k}(\sigma)=X_{k}^{*}(\sigma)$, $P(\sigma)^{*}
=P(\sigma)$, and subject to the relations
$$
P(\sigma )P(\tau)=P(\sigma \cap \tau), \;\;
P(\emptyset )=1
$$
$$
P(\sigma)X_{k}(\tau)=X_{k}(\tau)P(\sigma) \;\; {\rm iff} \;\; k\in \sigma
$$
i.e. $P(\sigma)$'s are projections which ``partially commute'' 
with the variables $X_{k}(\sigma)$. 
When equipped with the coproduct and the counit
$$
\widehat{\Delta} (X_{k}(\sigma))=X_{k}(\sigma)\otimes P(\sigma) + 
P(\sigma)\otimes X_{k}(\sigma)
$$
$$
\widehat{\Delta}(P(\sigma))=P(\sigma)\otimes P(\sigma),\;\;
\widehat{\epsilon}(X_{k}(\sigma))=0,\;\;
\widehat{\epsilon}(P(\sigma))=1,
$$
the algebra $\widehat{\cal B}$ becomes a unital *-bialgebra
called {\it filtered *-bialgebra}.
Therefore, we are in the position to study random walks [14]
and stochastic processes over *-bialgebras (see [1] and [18]).

We take a suitable state $\widehat{\phi}$ on $\widehat{\cal B}$, which
is obtained by lifting the tensor product state $\widetilde{\phi}^{\otimes 
\infty}$ on $\bigotimes _{k=1}^{\infty}{\bf C}\langle Y_{k}, P_{k}\rangle$ 
to $\widehat{\cal B}$ through the mapping which sends
each $X_{k}(\sigma)$ onto $Y_{k}$ and $P(\sigma)$ onto the tensor
product of $P_{k}$'s with $k\in \sigma$, where
$Y_{k}^{*}=Y_{k}$ and $P_{k}$ a projection.
The state $\widehat{\phi}$ is our noncommutative, ``filtered'' analog of
the classical product measure ($\sigma$'s play the role of 
filters due to ``partial commutations'').

The corresponding
convolution central limit theorem
(or discrete random walk),
which plays the role of a noncommutative
analog of the classical multivariate central limit theorem,
gives, under the usual normalization, 
pointwise convergence of the $N$-th convolution power
$$
\widehat{\phi}^{\star N}=\widehat{\phi}^{\otimes N}\circ \widehat{\Delta}^{N-1}
$$ 
where $\widehat{\Delta}^{N-1}$ is the $N-1$-th iteration of the 
coproduct $\widehat{\Delta}$.

The summands produced by iterating the coproduct are
called {\it filtered random variables} and can be viewed as
quantum analogs of independent {\it random vectors}.
It is important to note that
by taking suitable linear combinations (strongly convergent series
on the GNS pre-Hilbert space) of filtered random variables
we obtain $m$-free (free) random variables.
Thus all three basic notions of quantum independence 
in the axiomatic theory (tensor, free and Boolean,
see [3],[20]) are covered by this scheme. 

By considering random walks with continuous time,
or stochastic processes over the filtered *-bialgebra, we obtain
in the limit the vacuum expectation state in the 
multiple symmetric Fock space $\Gamma({\cal H})$, where 
$$
{\cal H}=L^{2}({\bf R}^{+})\otimes {\cal G}
$$
and ${\cal G}$ is a separable Hilbert space called the multiplicity
space. The GNS representation leads to 
{\it filtered fundamental operators} which are the 
CCR fundamental operators on $\Gamma(L^{2}({\bf R}^{+}) {\cal G})$, 
mutliplied by projections $P^{(\sigma)}$, where $P^{(\sigma)}$ is 
the second quantization of
the canonical projection onto subspaces of 
$L^{2}({\bf R}^{+})\otimes {\cal G}$
built from the modes (called {\it colors}) which belong to the 
set $\sigma$. 

Fundamental operators associated with different 
notions of independence can be expressed in terms
of the filtered ones. In particular, one can define bounded extensions
to $\Gamma(L^{2}({\bf R}^{+}) {\cal G})$ of
$m$-free creation and annihilation operators introduced in [5] as
strongly convergent series of filtered creation and annihilation operators, 
respectively.
This formalism enables us not only to
embed the free (or, full) Fock space over $L^{2}({\bf R}^{+})$ 
in $\Gamma({\cal H})$, but also decompose $\Gamma({\cal H})$
into an orthogonal sum of subspaces which are isomorphic
to the free Fock space. 

The corresponding {\it filtered
stochastic calculus} is developed in [12]
and it is, in fact, a generalization of the 
Hudson-Parthasarathy calculus [6] (see also [16]) 
on multiple symmetric Fock
spaces [15] and includes a new version of the 
free calculus, originally developed
for the Cuntz algebra [8], as well as gives
calculi for the hierarchy of $m$-free Brownian motions 
introduced in [4]. In that context, see also [7] and [17].

In Section 2 we introduce the filtered *-bialgebra which
sets the framework for a unified approach to noncommutative
probability.
This leads to filtered random variables, which are introduced
in the more general setting of unital *-algebras in Section 3.
Their combinatorics and the recurrence relation for the
product state is given in Section 4. Convolution 
limit theorems are proved in Section 5.
In Section 6 we introduce the filtered fundamental operators.
These are used for the GNS 
construction of the limit of a sequence of random walks 
on the filtered *-bialegbra in Section 7.
In Section 8 we determine the combinatorics
of general filtered white noises.
In Section 9 we study in more detail extensions of
the $m$-free and free creation and annihilation operators 
to all of $\Gamma({\cal H})$. A free Fock space
decomposition of $\Gamma({\cal H})$ is established.

We denote all scalar products
by $\langle ., . \rangle$ and identify operators and their
ampliations if no confusion arises.
\\[10pt]
\myownsection
\begin{center}
{\sc 2. Filtered bialgebras and convolutions}
\end{center}
In this section we discuss the bialgebras in our construction 
and the associated convolution. For general background on this,
we refer the reader to [14] and [18].

For simplicity, consider first the unital *-algebra 
${\bf C}[X]$ of polynomials in the variable $X=X^{*}$ 
endowed with the coproduct 
$$
\Delta :\; {\bf C}[X]\rightarrow {\bf C}[X]\otimes {\bf C}[X]
$$ 
given by
\begin{equation}
\label{2.1}
\Delta (X)=X\otimes 1 + 1 \otimes X
\end{equation}
and the counit $\epsilon : {\bf C}[X]\rightarrow {\bf C}$
given by $\epsilon (X)=0$.
This coproduct leads to the {\it classical convolution} of measures
and thus classical convolution of states.

Namely, if $\phi, \psi$ are states on ${\bf C}[X]$ associated with measures
$\mu$, $\nu$ on the real line, i.e.
$$
\phi(X^{n})=\int_{{\bf R}}x^{n}d\mu (x), \;\;\;
\psi(X^{n})=\int_{{\bf R}}x^{n}d\nu(x),
$$
then the convolution of states 
$$
\phi \star_{c} \psi  = \phi \otimes \psi \circ \Delta 
$$
corresponds to 
the classical convolution of measures $\mu \star_{c} \nu$ in the sense that
$$
\phi \star_{c} \psi (X^{n})=m_{n}(\mu \star_{c} \nu),
$$
where $m_{n}(\mu \star_{c} \nu)$ is the $n$-th moment of 
the measure $\mu\star_{c}\nu$.

The coproduct is a convenient tool to 
produce independent random variables [14]. 
Namely, by applying succesive iterations of $\Delta$ to $X$, we obtain
$$
\Delta^{N-1}(X)= \sum_{k=1}^{N}j_{l,N}(X)
$$
where $\Delta^{N}:=({\rm id}\otimes \Delta^{N-1})\circ \Delta$
for $N>1$ with $\Delta^{1}=\Delta$, and the summands
$$
j_{l,N}(X)=1^{\otimes (l-1)}\otimes X \otimes 1^{\otimes (N-l)},\;\;
1\leq l \leq N,
$$
can be viewed as independent random variables with respect to the state
$\phi^{\otimes N}$.

The so-called {\it Boolean convolution} can be obtained by
considering 
the unital *-algebra of polynomials in two noncommuting self adjoint
variables
${\bf C}\langle X', P \rangle$, where $P$ is a projection,
with the coproduct
$$
\Delta:\; {\bf C}\langle X', P \rangle \rightarrow {\bf C}\langle X', P\rangle
\otimes {\bf C}\langle X' , P \rangle , 
$$
given by
\begin{equation}
\label{2.2}
\Delta (X')=X'\otimes P + P \otimes X', \;\;
\Delta(P)=P\otimes P, \;\;
\end{equation}
and the counit 
$\epsilon (X')=0$, $\epsilon (P)=1$. 
It can be shown that this coproduct gives the Boolean convolution of states and
thus the Boolean convolution of measures [23]. 
This follows from the hierarchy of freeness construction 
[11], but a direct proof will be presented below.\\
\indent{\par}
{\sc Definition 2.1.}
If $\phi$ is a state on ${\bf C}[Y]$, where $Y^{*}=Y$,
then its {\it Boolean extension}
is the state on ${\bf C}\langle Y, P \rangle $, where $P$ is a projection,
given by the linear extension of
\begin{equation}
\label{2.3}
\widetilde{\phi}(P^{\alpha}Y^{n_1}PY^{n_2}P\ldots Y^{n_k}P^{\beta})
=\phi(Y^{n_1})\phi(Y^{n_2})\ldots \phi(Y^{n_k})
\end{equation}
where $\alpha , \beta \in \{0,1\}$ and $n_{1}, \ldots , n_{k}\in {\bf N}$, with
$\widetilde{\phi}(P)=1$. If $\phi$ is a state on the unital *-algebra
${\cal A}$, then its Boolean extension $\widetilde{\phi}$ to
the free product $\tilde{\cal A}={\cal A}*{\bf C}[P]$ (units identified)
is defined in an analogous way.\\
\indent{\par}
The Boolean extension of a state is a state since it is obtained
as the Boolean product of the state $\phi$ on ${\bf C}[Y]$ and the
unital *-homomorphism $h$ on ${\bf C}[P]$ given by $h(P)=h(1)=1$.

In order to have a unified model for both convolutions it is now enough
to incorporate both coproducts (2.1)-(2.2) into one scheme. 
This is done as follows.
The unital *-algebra ${\cal B}={\bf C}\langle X,X',P\rangle$
where $X=X^{*}$, $X'=X'^{*}$ and $P$ is a projection, 
endowed with the coproduct 
$\Delta:\; {\cal B}\rightarrow {\cal B}\otimes {\cal B}$
and counit $\epsilon: {\cal B}\rightarrow {\bf C}$ given by 
$$
\Delta (X)= X\otimes 1 + 1 \otimes X, \;\;
\Delta (X')=X'\otimes P + P \otimes X'
$$
$$
\Delta (P)=P\otimes P, \;\;\; \epsilon (X)=\epsilon (X')=0, \;\;
\epsilon (P)=1
$$
(in other words, $X'$ is $P$-primitive and $P$ is group-like),
becomes a unital *-bialgebra.  Both classical and Boolean convolutions are
recovered from $({\cal B}, \Delta , \epsilon)$, as we show below.\\
\indent{\par} 
{\sc Proposition 2.2.}
{\it  Let $\eta: \; {\cal B}\rightarrow {\bf C}\langle Y,P
\rangle $, where $Y=Y^{*}$ and $P^{2}=P=P^{*}$,
be the linear and multiplicative extension of}
$$
\eta (X)=\eta (X')=Y, \;\; \eta (P)=P, \;\; \eta (1)=1
$$
{\it and, for states $\phi ,\psi$ on ${\bf C}[Y]$, 
let ${\phi}_{0}=\widetilde{\phi} \circ \eta$, 
${\psi}_{0}=
\widetilde{\psi}\circ \eta$ with the convolution}
\begin{equation}
\label{2.4}
{\phi}_{0}\star {\psi}_{0}=
{\phi}_{0}\otimes {\psi}_{0}\circ \Delta
\end{equation}
{\it where $\Delta $ is the coproduct for ${\cal B}$. Then the restrictions of
${\phi}_{0}\star {\psi}_{0}$ to ${\bf C}[X]$ and ${\bf C}[X']$,
respectively, agree with $\phi \star_{c} \psi$ and $\phi\star_{B}\psi$, 
respectively.}\\[5pt]
{\it Proof.}
First of all, note that ${\phi}_{0}$ and ${\psi}_{0}$ are 
states since $\eta $ is a unital *-homomorphism.
That the restriction of the convolution (2.4) to ${\bf C}[X]$ gives
classical convolution, is obvious. In turn, the statement concerning
the Boolean convolution follows from the fact that the subalgebras
${\bf C}[Y]\otimes P$, $P\otimes {\bf C}[Y]$ of ${\bf C}\langle Y, P \rangle
\otimes {\bf C}\langle Y , P \rangle $
are Boolean independent with respect to
the state $\widetilde{\phi}\otimes \widetilde{\psi}$. This fact
can be easily seen from the following calculation:
$$
\widetilde{\phi}\otimes \widetilde{\psi}(Y^{k_{1}}\otimes P)(P\otimes
Y^{k_{2}})(Y^{k_{2}}\otimes P)(P\otimes Y^{n_{2}}) \ldots)
$$
\begin{eqnarray*}
&=& \widehat{\phi}(Y^{k_1}PY^{k_2}P \ldots )
\widehat{\psi}(PY^{n_1}PY^{n_2}\ldots )\\
&=&\phi(Y^{k_{1}})\psi(Y^{n_{1}})\phi(Y^{k_{2}})\psi (Y^{n_{2}})\ldots
\end{eqnarray*}
A more general setting of the Bollean product of states was given in [11].
\hfill $\Box$\\
\indent{\par}
The quadruple $({\cal B}, \Delta , \epsilon , \phi_{0})$ can be called 
the {\it random walk} (we follow Majid [14] in this terminology) 
on the {\it pair of quantum planes}.
Note that on the quantum probability space level
we may also study the pair $({\bf C}\langle Y, P \rangle , \widetilde{\phi})$, 
which corresponds to (polynomial functions on) a 
{\it pair of quantum real lines}.

Let us consider now the multivariate generalization of the *-bialgebra
${\cal B}$.
In classical probability, the multivariate case in an algebraic formulation
would be reached if we considered the unital *-algebra 
${\bf C}[X_{k}; k\in {\bf N}]$ of polynomials in commuting variables
$(X_{k})_{k\in {\bf N}}$, with the classical coproduct 
\begin{equation}
\label{2.5}
\Delta (X_{k})=X_{k}\otimes 1 + 1 \otimes X_{k}
\end{equation}
and the counit $\epsilon (X_{k})=0$. 

Let us now define a quantum analog of this multivariate *-bialgebra, of which
the bialgebra ${\cal B}$ is the ``one-dimensional'' version.
Thus, introduce the unital *-algebra
$$
\widehat{\cal B}=
{\bf C}\langle X_{k}(\sigma), P(\sigma);\; 
k\in {\bf N}, \sigma\in {\cal P}({\bf N}) \rangle/J
$$
where $X_{k}(\sigma)=X_{k}^{*}(\sigma)$, $P(\sigma)^{*}=P(\sigma)$,
${\cal P}({\bf N})$ is the power set of ${\bf N}$ 
and $J$ is the two-sided 
ideal generated by the relations
\begin{equation}
\label{2.6}
P(\sigma )P(\sigma ')=P(\sigma \cap \sigma '), \;\;
P(\emptyset)=1
\end{equation}
\begin{equation}
\label{2.7}
P(\sigma)X_{k}(\tau)=X_{k}(\tau)P(\sigma) \;\; {\rm iff}\;\; k\in \sigma
\end{equation}
i.e. the projection associated with the set $\sigma$ 
(we call $\sigma$ a { \it filter})
``filters through'' the variables
$X_{k}(\sigma )$ if the index $k\in \sigma$.\\
\indent{\par}
{\sc Proposition 2.3.}
{\it The algebra $\widehat{\cal B}$, equipped with the coproduct
$\widehat{\Delta}: 
\widehat{\cal B}\rightarrow \widehat{\cal B}\otimes \widehat{\cal B}$
and the counit $\widehat{\epsilon}: \widehat{\cal B}\rightarrow {\bf C}$
given by}
\begin{equation}
\label{2.8}
\widehat{\Delta} (X_{k}(\sigma))=X_{k}(\sigma)\otimes P(\sigma) + 
P(\sigma)\otimes X_{k}(\sigma) 
\end{equation}
\begin{equation}
\label{2.9}
\widehat{\Delta}(P(\sigma))
=P(\sigma)\otimes P(\sigma), \;\;
\widehat{\epsilon} (X_{k}(\sigma))=0, \;\;
\widehat{\epsilon} (P(\sigma))=1
\end{equation}
{\it for all $k$ and $\sigma$, is a unital *-bialgebra called
filtered *-bialgebra.}\\[5pt]
{\it Proof.} The coproduct and the counit preserve the relations
(2.6)-(2.7). \hfill $\Box$\\
\indent{\par}
By iterating this coproduct, called {\it filtered coproduct}, we obtain 
the sum
\begin{equation}
\label{2.10}
\widehat{\Delta}^{N-1}(X_{k}(\sigma))=
\sum_{l=1}^{N}\widehat{j}_{l,N}(X_{k}(\sigma))
\end{equation}
of $P(\sigma)$-deformed random variables
\begin{equation}
\label{2.11}
\widehat{j}_{l,N}(X_{k}(\sigma))=P(\sigma)^{\otimes l-1}\otimes X_{k}(\sigma) 
\otimes P(\sigma)^{\otimes (N-l)},
\end{equation}
where $\sigma\in {\cal P}({\bf N})$, $k\in {\bf N}$, $1\leq l \leq N$,
$N\in {\bf N}$.\\
\indent{\par}
{\it Remark.}
The power set ${\cal P}({\bf N})$ can be put in one-to-one 
correspondence
with the Cantor set $CS$ since any infinite sequence 
$(q_{1}, q_{2}, q_{3}, \ldots )$ of $0$'s and $2$'s
gives a unique number $q\in CS$ with the ternary expansion given by
$$
q=\sum_{n=0}^{\infty}\frac{q_{n}}{3^{n}}
$$
which allows us to identify the set $\sigma$ with a number $q\in CS$
according to the rule
$$
n\in \sigma \;\; {\rm iff} \;\; q_{n}=0
$$
Thus, the index set which gives our discrete quantum 
deformation of the coproduct 
(2.5) plays a role similar to the interval 
$[0,1]$ in the case of $q$-deformed quantum groups like $SU_{q}(2)$ or
$U_{q}(su(2))$ (see [14],[10],[13]). \\
\indent{\par}
When we go over from quantum groups to quantum probability spaces, 
we identify $X_{k}(\sigma)$'s for all different $\sigma$ and fixed $k$, 
as we did by using the map $\eta$ in the ``one-dimensional''case
of Proposition 2.2. 
This is done in order to include different notions of independence in one
scheme. For that purpose, we consider the mapping
\begin{equation}
\label{2.12}
\widehat{\eta}: \; \widehat{\cal B}\rightarrow \bigotimes_{k=1}^{\infty}
{\bf C}\langle Y_{k}, P_{k} \rangle
\end{equation}
given by the linear and multiplicative extension of
\begin{eqnarray}
\label{2.13}
\widehat{\eta} (P(\sigma))&=& P_{1}^{q_1}\otimes P_{2}^{q_2}\otimes
\ldots \otimes P_{k}^{q_k}\otimes \ldots\\
\label{2.14}
\widehat{\eta} (X_{k}(\sigma))&=&1_{1}\otimes 1_{2}\otimes \ldots \otimes 1_{k-1}\otimes 
Y_{k}\otimes 1_{k+1}\otimes \ldots 
\end{eqnarray}
where ${\bf C}\langle Y_{k}, P_{k} \rangle$ is the $k$-th copy of 
${\bf C}\langle Y, P\rangle $ and the sequence $(q_{1},q_{2},q_{3},\ldots )$
represents $\sigma$. 
The infinite tensor product is taken with respect to the set 
$\{1_{k},P_{k}, k\in {\bf N}\}$ (see [5] for the formal definition).

It can be seen that $\eta$ is a unital *-homomorphism.
Therefore, for given state $\phi$ on ${\bf C}[Y]$, 
the functional 
\begin{equation}
\label{2.15}
\widehat{\phi}=\widetilde{\phi}^{\otimes \infty}\circ \widehat{\eta},
\end{equation}
is a state on $\widehat{\cal B}$. It plays the role of a 
noncommutative analog of a vector state in classical probability.
A generalization to vector states corresponding to products of 
different measures is immediate. It is enough to take
\begin{equation}
\label{2.16}
\widehat{\phi}=\bigotimes_{k=1}^{\infty}\widetilde{\phi}_{k}
\circ \widehat{\eta},
\end{equation}
where $\phi_{k}$ is a state on ${\bf C}[Y_{k}]$, $k\in {\bf N}$.

When we take (2.15) (or, (2.16)),
the quadruple $(\widehat{\cal B},\widehat{\Delta} , 
\widehat{\epsilon}, \widehat{\phi})$
will give our stationary 
{\it filtered random walk}, using the terminology of Majid [14], 
which carries two structures, that of the unital * -bialgebra and 
that of the quantum probability space.
The corresponding convolution of states
\begin{equation}
\label{2.17}
\widehat{\phi}\star \widehat{\psi}=\widehat{\phi}\otimes \widehat{\psi}
\circ \widehat{\Delta}
\end{equation}
(or, products of states)
will be called the {\it filtered convolution}.

One of the main motivations to study the filtered *-bialgebras, 
convolutions, random walks and stochastic processes comes from the 
following result, which follows from our previous work [11], although
it has not been stated there in terms of the filtered convolution.\\
\indent{\par}
{\sc Proposition 2.4.}
{\it If $\widehat{\phi}$ and $\widehat{\psi}$ are 
of the form (2.16) and $i$ is the unital *-homomorphism}
$$
i:\; {\bf C}\langle X_{k}, k\in {\bf N}\rangle\rightarrow \widehat{\cal B},\;\;\;
i(X_{k})=X_{k}({\bf N}),
$$
{\it then $\widehat{\phi}\star \widehat{\psi}\circ i$ 
agrees with the classical convolution of products of states. 
If $\widehat{\phi}$ and $\widehat{\psi}$
are of the form (2.15) and $i^{(m)}$ denotes the unital *-homomorphism}
$$
i^{(m)}:\;\; {\bf C}[X]\rightarrow \widehat{\cal B},\;\;
i^{(m)}(X)=\sum_{k=1}^{m}(X_{k}(k)-X_{k}(k-1))
$$
{\it where $X_{k}(p)=X_{k}(\{1, \ldots , p-1\})$, 
then $\widehat{\phi}\star \widehat{\psi}\circ i^{(m)}$ 
agrees with the additive $m$-free convolution of states 
$\phi\star_{m}\psi$ on ${\bf C}[X]$, $1\leq m \leq \infty$.}
\\[5pt]
{\it Proof.}
The statement concerning the classical convolution of products of states
is obvious since $P({\bf N})=1$. In turn, the second part of the 
proposition is non-trivial and 
follows from the construction of $m$-free product states
and the associated *-bialgebras (see [11], Section 5, where we also refer 
the reader for the definition of the $m$-free convolution). 
\hfill $\Box$ \\
\indent{\par}
Since
${\bf C}\langle Y, P\rangle $ can be viewed as a 
{\it quantum pair of real lines}, on the 
quantum probability space level we can interpret our object of 
interest as (polynomial functions on) the 
product of infinitely many quantum pairs of real lines,
which is our non-commutative analog of ${\bf R}^{\infty}$.
On the bialgebra level, we have a bigger object
since every variable $X_{k}$ admits a family of 
different convolutions.\\[10pt]
\myownsection
\begin{center}
{\sc 3. Filtered random variables}
\end{center}
In this section we introduce filtered random variables 
which are our noncommutative analogs of independent random vectors
in the general setting of arbitrary unital *-algebras. 

In analogy to the classical case, we obtain them
by iterating the coproduct $\widehat{\Delta}$ as in (2.10). 
Then, we embed $\widehat{j}_{l,N}(X_{k}(\sigma)$
into $\widehat{\cal B}^{\infty}$ to get
\begin{equation}
\label{3.1}
P(\sigma)^{\otimes (l-1)}\otimes X_{k}(\sigma) \otimes P(\sigma)^{\otimes \infty}
\end{equation}
and generalize these to the arbitrary unital *-algebras.
In order to do that, write (3.1) as the product
$$
(1^{\otimes (l-1)}\otimes X_{k}(\sigma) \otimes 1^{\otimes \infty})
(P(\sigma)^{\otimes (l-1)}\otimes 1 \otimes P(\sigma)^{\otimes \infty})
$$
of an ampliation of
$X_{k}(\sigma)$ into $\widehat{\cal B}^{\otimes \infty}$ and a projection
indexed by $\sigma$. This shows that the definitions given below are
a natural generalization of those of Section 2.

Let $({\cal A}_{l})_{l\in L}$ be a family of unital
*-algebras with units $1_{l}$ and let 
$(\phi_{l})_{l\in L}$ be the corresponding family of states. 
Consider a noncommutative probability space
$(\widehat{\cal A}_{1}, \widehat{\Phi}_{1})$, where
$$
\widehat{\cal A}_{1}
= 
\bigotimes_{l\in L}\widetilde{\cal A}_{l}^{\otimes \infty},\;\;\;
\widehat{\Phi}_{1}
=
\bigotimes_{l\in L}\widetilde{\phi}_{l}^{\otimes \infty},
$$
and $\widetilde{\cal A}_{l}={\cal A}_{l}*{\bf C}[P_{l}]$ is the free
product with identified units, $P_{l}$ being a projection,
whereas $\widetilde{\phi}$ is the Boolean extension of $\phi$
(Definition 2.1).
The infinite tensor products are understood as in 
[5], with the canonical involution. 
This noncommutative probability space will be called the 
{\it multiple probability space} associated with the considered
family of probability spaces
since each of them appears 
infinitely many times in the considered tensor products.
We will refer to those copies as {\it colors}. Roughly speaking,
$\widetilde{\cal A}_{l}^{\otimes \infty}$ 
and $\widetilde{\phi}_{l}^{\otimes \infty}$ correspond to
$\bigotimes_{k=1}^{\infty}{\bf C}\langle Y_{k} , P_{k} \rangle$ and
$\widetilde{\phi}^{\otimes \infty}$ for each $l\in L$, respectively
of Section 2.

If $({\cal H}_{l}, \pi_{l} , \Omega_{l})$ 
is the GNS triple for the pair $({\cal A}_{l}, \phi_{l})$, 
then $({\cal H}_{l}, \widetilde{\pi}_{l} , \Omega_{l})$ is the
GNS triple for $(\widetilde{\cal A}_{l}, \widetilde\phi_{l})$, $l\in L$,
where $\widetilde{\pi}_{l}$ agrees $\pi_{l}$ on
${\cal A}_{l}$ and $\widetilde{\pi}_{l}(P_{l})$ is 
the projection onto the cyclic vector $\Omega_{l}$.
For convenience, we can identify 
$x\in {\cal A}_{l}$ with $\pi_{l}(x)$, 
$P_{l}$ with $P_{\Omega_{l}}$ and $\phi_{l}$ with the expectation state 
$\langle \Omega_{l}, . \Omega_{l}\rangle$ (see [5]).

Guided by (3.1), from 
projections $P_{m}$ we construct projections ${\bf P}(l,\sigma)$ 
to be elementary tensors in $\widehat{\cal A}_{1}$ with components 
$$
{\bf P}(l,\sigma)_{m,k}=
\left\{ 
\begin{array}{cc}
P_{m} & {\rm if} \;\; m\neq l\;\; {\rm and}\;\; k\notin \sigma\\
1_{m} & {\rm otherwise}
\end{array}
\right.
$$
where $\sigma\in {\cal P}({\bf N})$ and $l\in L$. 
In the case when 
$\sigma=\{1, \ldots , r-1\}$, 
we will write ${\bf P}(l,r)={\bf P}(l,\sigma)$.\\
\indent{\par}
{\sc Definition 3.1.}
By {\it filtered random variables} we will understand
elements of $\widehat{\cal A}_{1}$ which are of the form
\begin{equation}
\label{3.2}
{\bf XP}
\end{equation}
where ${\bf X}={\bf X}(l,k)$ is the 
$(l,k)$-th ampliation of $x\in {\cal A}_{l}$ into $\widehat{\cal A}_{1}$
and ${\bf P}={\bf P}(l,\sigma)$, where
$l\in L$, $k\in {\bf N}$, $\sigma\in {\cal P}({\bf N})$.
In particular, the unit 
${\bf 1}=\bigotimes_{l\in L}1_{l}^{\otimes \infty}$
is a filtered random variable.\\
\indent{\par}
If $L={\bf N}$, filtered random variables can be represented 
as infinite matrices. 
Assume that $k$ numbers the rows and $l$ numbers the columns.
For instance, for $x\in {\cal A}_{2}$, ${\bf X}={\bf X}(2,3)$ 
and ${\bf P}={\bf P}(2,4)$
we have
$$
{\bf XP}=
\left(
\begin{array}{llll}
1_{1} & 1_{2} & 1_{3} & \ldots \\
1_{1} & 1_{2} & 1_{3} & \ldots \\
1_{1} & x & 1_{3} & \ldots \\
P_{1} & 1_{2} & P_{3} & \ldots \\
P_{1} & 1_{2} & P_{3} & \ldots \\
. & . & . & \ldots 
\end{array}
\right)
$$
Note that if $L={\bf N}$, multiplication of filtered random variables
corresponds to Schur's multiplication of matrices. \\
\indent{\par}
{\sc Definition 3.2.}
Let $\widehat{\cal A}$ be the unital *-subalgebra 
of $\widehat{\cal A}_{1}$ generated by
all filtered random variables and let 
$\widehat{\Phi}=\widehat{\Phi}_{1}|\widehat{\cal A}$.
The noncommutative probability space 
$(\widehat{\cal A}, \widehat{\Phi})$ will be called
the {\it filtered probability space} associated with
$({\cal A}_{l}, \phi_{l})_{l\in L}$ and the state
$\widehat{\Phi}$ will be called the {\it filtered
product} of $(\phi_{l})_{l\in L}$. 
The unital *-subalgebras of $\widehat{\cal A}$,
$$
\widehat{\cal A}_{l}=
\langle {\bf XP}|\; {\bf X}={\bf X}(l,k),\; 
{\bf P}={\bf P}(l,\sigma),\;x\in {\cal A}_{l},
\; k\in {\bf N}, \; 
\sigma\in {\cal P}({\bf N}) \rangle ,\;\;l\in L
$$
will be called {\it filtered with respect to $\widehat{\Phi}$}.\\
\indent{\par}
{\it Example 1.}
Let $*_{l\in L}{\cal A}_{l}$ denote the free product of 
$({\cal A}_{l})_{l\in L}$ with non-identified units.
Fix $k\in {\bf N}$, $\sigma\in {\cal P}({\bf N})$ and define a 
*-homomorphism
$$
j^{(k,\sigma)}:\;\; *_{l\in L}{\cal A}_{l}\rightarrow \widehat{\cal A}
$$
as the linear extension of
$$
j^{(k,\sigma)}(x_{1}\ldots x_{n})
=
{\bf X}_{1}(l_{1},k)
{\bf P}(l_{1},\sigma)\ldots 
{\bf X}_{n}(l_{n},k)
{\bf P}(l_{n},\sigma)
$$
for $x_{i}\in {\cal A}_{l_{i}}$, $l_{1}\neq l_{2} \neq \ldots \neq l_{n}$.
Let us define the mapping
$$
i:\;\; *_{l\in L}{\cal A}_{l}\rightarrow \bigotimes_{l\in L}{\cal A}_{l} 
$$
as the linear extension of
$$
i(x_{1}\ldots x_{n})=i_{l_{1}}(x_{1}) \ldots i_{l_{n}}(x_{n})
$$
where $i_{l}$ are canonical *-homomorphic embeddings of ${\cal A}_{l}$
into $\bigotimes_{l\in L}{\cal A}_{l}$. Then, 
\begin{eqnarray*}
\widehat{\Phi}\circ j^{(k,\sigma)}
=\left\{
\begin{array}{ccc}
\bigotimes_{l\in L}\phi_{l}\circ i &{\rm if} &k\in \sigma\\[5pt]
*^{B}_{l\in L}\phi_{l} & {\rm if} &k\notin \sigma
\end{array}
\right.
\end{eqnarray*}
where $*^{B}_{l\in L}\phi_{l}$ denotes the Boolean (or, $1$-free)
product of states $(\phi_{l})_{l\in L}$ on $*_{l\in L}{\cal A}_{l}$.
In other words, for fixed $k,\sigma$, the set 
$\{{\bf X}(l,k){\bf P}(l,\sigma): l\in L\}$ is a family
tensor independent r.v. if 
$k\in \sigma$ and Boolean independent r.v.
if $k\notin \sigma$ (see [7]). \\
\indent{\par}
{\it Example 2.}
As we showed in [11], 
the Boolean product is just the first-order approximation of the 
free product of states in free probability [24].
Higher order approximations given by the modified hierarchy
of $m$-free products [5] (the case with non-identified units of
the usual hierarchy of freeness of [11])
can also be obtained from the filtered product.
Namely, let $m\in {\bf N}$, and define
$$
\bar{j}^{(m)}:\;\; *_{l\in L}{\cal A}_{l} \rightarrow \widehat{\cal A}
$$
as the linear extension of
$$
\bar{j}^{(m)}(x_{1}\ldots x_{n})=
\bar{j}_{l_{1}}^{(m)}(x_{1})\ldots 
\bar{j}_{l_{n}}^{(m)}(x_{n})
$$
where $x_{i}\in {\cal A}_{l_{i}}$, $l_{1}\neq l_{2} \neq \ldots \neq l_{n}$,
and 
\begin{equation}
\label{3.3}
\bar{j}_{l}^{(m)}(x)=\sum_{k=1}^{m}{\bf X}(l,k)({\bf P}(l,k)-{\bf P}(l,k-1))
\end{equation}
for $x\in {\cal A}_{l}$. Then
$$
\widehat{\Phi}\circ \bar{j}^{(m)}=*_{l\in L}^{(m)}\phi_{l}
$$
where $*_{l\in L}^{(m)}\phi_{l}$ denotes the modified $m$-free product of 
states. 
If $m=\infty$, the series given by the representation of 
(3.3) converges strongly on the GNS pre-Hilbert space (see [5]).
Moreover,
$$
\bar{j}_{l}^{(\infty)}({\bf 1}_{l})={\bf 1},\;\; l\in L
$$
and thus $\widehat{\Phi}\circ \bar{j}^{(\infty)}$ is well-defined
on the free product of ${\cal A}_{l}$, $l\in L$, with identified units
and agrees on it with the free product of states (for details, 
see [5]). Thus, the variables
$\bar{j}_{l}^{(m)}(x)$, $x\in {\cal A}_{l}$, $l\in L$, are $m$-free
random variables for $m\in {\bf N}$ and 
free random variables if $m=\infty$.\\[10pt]
\myownsection
\begin{center}
{\sc 4. Combinatorics}
\end{center}
Let us now introduce a new class of partitions which is crucial to the 
combinatorics of filtered random variables.\\
\indent{\par}
{\sc Definition 4.1.}
Let $\vec{k}=(k_{1}, \ldots , k_{n})$ and 
$\vec{\sigma}=(\sigma_{1}, \ldots , \sigma_{n})$
be color and filter tuples of natural numbers and sets of
natural numbers, respectively.
A partition $R=\{R_{1}, \ldots , R_{q}\}$ 
of the set $\{1, \ldots , n\}$ will be called 
$(\vec{k},\vec{\sigma})$-{\it adapted}
if and only if it satisfies the conditions\\
\indent{\par} (A1). 
$
\forall \;1\leq q \leq n\;
\forall \;i,j \in R_{q}\;\; \;
k_{i} =k_{j}$ \\
\indent{\par} (A2). 
${\rm If} \;\;i<m<j,\;\;{\rm where}\;\;
i,j\in R_{q}\;\;{\rm and}\;\; m\notin R_{q}, 
{\rm then}\;\; k_{i}=k_{j}\in \sigma_{m}$.\\[10pt]
The collection of all $(\vec{k},\vec{\sigma})$-adapted partitions 
(pair partitions) will be denoted by ${\cal P}_{n}(\vec{k},\vec{\sigma})$
(${\cal P}_{n}^{\rm pair}(\vec{k},\vec{\sigma})$).
The partitions of $\{1, \ldots , n\}$ which are
not $(\vec{k},\vec{\sigma})$-adapted will be called 
$(\vec{k},\vec{\sigma})$-{\it non-adapted}. \\
\indent{\par}
In other words, ${\cal P}_{n}(\vec{k},\vec{\sigma})$ is the subset
of all partitions ${\cal P}_{n}$ 
of $\{1, \ldots , n\}$ which are adapted to
the tuples $\vec{k}$ and $\vec{\sigma}$ in the following sense: 
(A1) colors corresponding to the elements of 
the same block have to match, 
(A2) between the elements of a given block there are no
filters associated with other blocks which separate them.  
In particular, if $k_{i}=k$, $\sigma_{i}=\sigma$ for all $1\leq i \leq n$, 
then the two extreme cases are given by 
$$
{\cal P}_{n}(\vec{k},\vec{\sigma})=
\left\{
\begin{array}{lll}
{\cal P}_{n} &{\rm if}& k\in \sigma\\
{\cal P}_{n}^{\rm int}& {\rm if} & k\notin \sigma
\end{array}
\right.
$$
where ${\cal P}_{n}^{\rm int}$ denotes the interval partitions
of $\{1, \ldots , n\}$. In turn, if $\sigma_{i}={\bf N}$ for all
$i=1, \ldots , n$, then 
$$
{\cal P}_{n}(\vec{k}, \vec{\sigma})=
{\cal P}_{n}(\vec{k})
$$
where ${\cal P}_{n}(\vec{k})$ denotes all partitions $R$
of the set $\{1, \ldots , n\}$
such that $k_{i}=k_{j}$ iff $ i,j$ belong to the same block of $R$
(this corresponds to the classical multivariate case).\\
\indent{\par}
{\sc Definition 4.2.}
If $R$ is $(\vec{k},\vec{\sigma})$-non-adapted, then 
the unique coarsest subpartition of $R$ which is 
$(\vec{k},\vec{\sigma})$-adapted will be denoted by $R(\vec{k},\vec{\sigma})$.
\\
\indent{\par}
{\it Examples.}
Consider the partition
$$
R=\{\{1,3,5\}, \{2,4\}\}
$$ 
of $\{1, \ldots , 5\}$ and let the color tuple
be given by $\vec{k}=(1,1,2,1,1)$. 
Then $R$ is not $(\vec{k}, \vec{\sigma})$-- adapted
for any $\vec{\sigma}$ since it does not satisfy (A1).
If we take now
the filter tuple $\vec{\sigma}=(\sigma_{1}, \ldots , \sigma_{5})$ given by 
$\sigma_{i}=\{1, \ldots , r_{i}-1\}$, with $r_{1}=r_{3}=r_{5}=1$,
$r_{2}=r_{4}=2$, then 
$$
R(\vec{k},\vec{\sigma})=\{\{1\}, \{2\} ,\{3\}, \{4\}, \{5\}\}.
$$
In turn, if we take
$\vec{\tau}=(\tau_{1},\ldots, \tau_{5})$, where
$\tau_{i}=\{1, \ldots , s_{i}-1\}$ and
$s_{1}=s_{5}=1$ and $s_{2}=s_{3}=s_{4}=2$, then
$$
R(\vec{k},\vec{\tau})=\{\{1, 5\}, \{2,4\}, \{3\}\}.
$$
\indent{\par}
We will see that partitions which are not $(\vec{k}, \vec{\sigma})$-- adapted
are less important since they do not survive in the limit theorems.
Therefore, there is an analogy with the crossing and non-crossing
partitions in free probability, the $(\vec{k},\vec{\sigma})$-- adapted
playing a similar role to non-crossing partitions, whereas the non--
$(\vec{k}, \vec{\sigma})$--adapted behave like crossing partitions.

Let us give a recurrence formula for moments of filtered random
variables, or ``filtered moments''.
It is convenient to introduce the following notions.\\
\indent{\par}
{\sc Definition 4.3.}
Given a tuple of pairs $((l_{1},k_{1}), \ldots (l_{n},k_{n}))$, 
we will say that 
$(l_j, k_j)$ is a {\it singleton} if $(l_{j},k_{j})\neq (l_{i},k_{i})$
for all $i\neq j$. If $(l_{i},k_{i})=(l_{j},k_{j})$ for $i<j$ such that
there is no $i<r<j$ for which $(l_{r},k_{r})=(l_{i},k_{i})$ 
and there exists
$i<m<j$ such that $l_{m}\neq l_{i}$ and $k_{i}\notin \sigma_{m}$,
then we will say that the filter $\sigma_{m}$ {\it separates}
$(l_i,k_i)$ and $(l_j,k_j)$.\\
\indent{\par} 
{\sc Proposition 4.4.}
{\it Let ${\bf X}_{i}:={\bf X}(l_{i},k_{i})$, ${\bf P}_{i}:=
{\bf P}(l_{i},\sigma_{i})$, 
$1\leq i \leq n$, where
$x_{i}\in {\cal A}_{l_{i}}$, $l_{i}\in L$, $k_{i}\in {\bf N}$,
$\sigma_{i}\in {\cal P}({\bf N})$, with $n\in {\bf N}$,
and $(l_{1},k_{1})\neq (l_{2},k_{2})\neq \ldots \neq (l_{n},k_{n})$.
Then}
$$
\widehat{\Phi}({\bf X}_{1}{\bf P}_{1}{\bf X}_{2}{\bf P}_{2}\ldots 
{\bf X}_{n}{\bf P}_{n})=
\widehat{\Phi}({\bf X}_{1}{\bf P}_{1})
\widehat{\Phi}({\bf X}_{2}{\bf P}_{2}\ldots {\bf X}_{n}{\bf P}_{n})
$$
{\it if $(l_{1},k_{1})$ is a singleton, or if there exists 
a filter $\sigma_{m}$ which separates $(l_{1},k_{1})$ and 
$(l_{r},k_{r})$, where
$r$ is the first index for which $(l_{1},k_{1})=(l_{r},k_{r})$, 
and otherwise}
$$
\widehat{\Phi}({\bf X}_{1}{\bf P}_{1}{\bf X}_{2}{\bf P}_{2}\ldots 
{\bf X}_{n}{\bf P}_{n})=
\widehat{\Phi}({\bf X}_{2}{\bf P}_{2}\ldots {\bf X}_{1}{\bf X}_{r}
{\bf P}_{r}\ldots {\bf X}_{n}{\bf P}_{n})
$$
{\it Proof.}
These formulas follow from the definition of filtered random variables.
\hfill $\Box$\\
\indent{\par}
{\sc Proposition 4.5.}
{\it Under the assumptions of Proposition 4.4,}
\begin{eqnarray*}
\widehat{\Phi}
(
{\bf X}_{1}{\bf P}_{1}
\ldots 
{\bf X}_{n}{\bf P}_{n}
)
&=&
\widehat{\Phi}({\bf X}_{B_{1}})\ldots \widehat{\Phi}({\bf X}_{B_{r}})\\
&=& 
\phi_{l(B_{1})}(x_{B_{1}})\ldots \phi_{l(B_{r})}(x_{B_{r}})
\end{eqnarray*}
{\it where $R$ is the partition associated with the tuple 
$(l_{1}, \ldots , l_{n})$,
$B_{1}, \ldots B_{r}$ are the blocks of 
$R(\vec{k}, \vec{\sigma})$, ${\bf X}_{B}=\prod_{j\in B}{\bf X}_{j}$ and
$x_{B}=\prod_{j\in B}x_{j}$ are products taken in the natural order, and
$l(B)$ is the index $l\in L$ associated with block $B$.}\\[5pt] 
{\it Proof.}
This is a straightforward consequence of Proposition 4.4 and the fact
that if $j_{1}, \ldots , j_{r}$ are elements of the same block 
$B\in R(\vec{k},\vec{\sigma})$, then 
$$
\widehat{\Phi}({\bf X}_{j_{1}}{\bf P}_{j_{1}}\ldots 
{\bf X}_{j_{r}}{\bf P}_{j_{r}})=
\widehat{\Phi}({\bf X}_{j_{1}}\ldots {\bf X}_{j_{r}}) 
=\phi_{l(B)}(x_{j_{1}}\ldots x_{j_{r}})
$$
(Definition 4.2 is crucial here). \hfill $\Box$\\
\indent{\par}
{\it Example.}
Let $L={\bf N}$ and take $x_{1},x_{3}\in {\cal A}_{1}$
and $x_{2}, x_{4} \in {\cal A}_{2}$. Then 
the partition $R=\{ \{1,3\}, \{2,4\}\}$ is 
associated with the tuple $(l_{1},l_{2},l_{3},l_{4})=(1,2,1,2)$. 
Let us consider two cases of color and filter tuples:
(i) $\vec{k}=(1,1,1,1)$, $\vec{\sigma}=(1,2,2,1)$ and 
(ii) $\vec{m}=(1,1,1,1)$, $\vec{\tau}=(1,2,1,1)$. Then the 
corresponding ``filtered moments'' can be obtained by refinement of
$R$ and represented in terms of diagrams. By $p^{k,r}$ we 
understand the number $p$ associated with color $k$ and filter
$\{1, \ldots , r-1\}$. \\
\unitlength=1mm
\special{em:linewidth 0.4pt}
\linethickness{0.4pt}
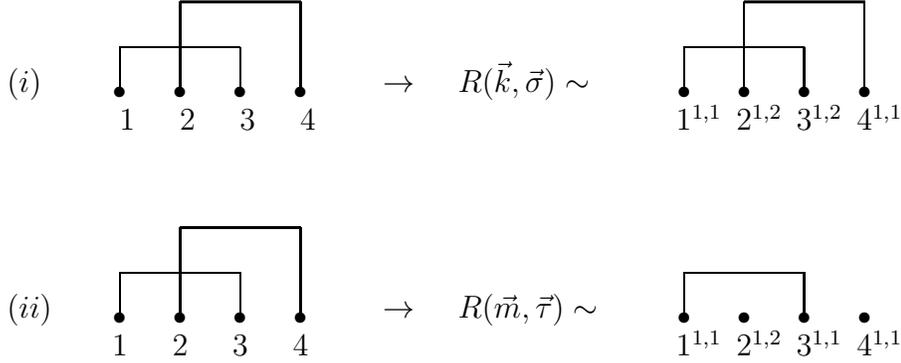
\begin{picture}(120.00,60.00)(-15.00,5.00)
\put(18.00,56.00){\line(1,0){16.00}}
\put(10.00,50.00){\line(1,0){16.00}}
\put(18.00,56.00){\line(0,-1){12.00}}
\put(34.00,56.00){\line(0,-1){12.00}}
\put(-5.00,44.00){$(i)$}
\put(10.00,50.00){\line(0,-1){6.00}}
\put(26.00,50.00){\line(0,-1){6.00}}
\put(10.00,44.00){\circle*{1.50}}
\put(18.00,44.00){\circle*{1.50}}
\put(26.00,44.00){\circle*{1.50}}
\put(34.00,44.00){\circle*{1.50}}
\put(10.00,39.00){$1$}
\put(18.00,39.00){$2$}
\put(26.00,39.00){$3$}
\put(34.00,39.00){$4$}
\put(45.00,44.00){$\rightarrow$}
\put(55.00,44.00){$R(\vec{k},\vec{\sigma})\sim$}
\put(93.00,56.00){\line(1,0){16.00}}
\put(85.00,50.00){\line(1,0){16.00}}
\put(93.00,56.00){\line(0,-1){12.00}}
\put(109.00,56.00){\line(0,-1){12.00}}
\put(85.00,50.00){\line(0,-1){6.00}}
\put(101.00,50.00){\line(0,-1){6.00}}
\put(85.00,44.00){\circle*{1.50}}
\put(93.00,44.00){\circle*{1.50}}
\put(101.00,44.00){\circle*{1.50}}
\put(109.00,44.00){\circle*{1.50}}
\put(84.00,39.00){$1^{1,1}$}
\put(92.00,39.00){$2^{1,2}$}
\put(100.00,39.00){$3^{1,2}$}
\put(108.00,39.00){$4^{1,1}$}
\put(18.00,26.00){\line(1,0){16.00}}
\put(10.00,20.00){\line(1,0){16.00}}
\put(18.00,26.00){\line(0,-1){12.00}}
\put(34.00,26.00){\line(0,-1){12.00}}
\put(-5.00,14.00){$(ii)$}
\put(10.00,20.00){\line(0,-1){6.00}}
\put(26.00,20.00){\line(0,-1){6.00}}
\put(10.00,14.00){\circle*{1.50}}
\put(18.00,14.00){\circle*{1.50}}
\put(26.00,14.00){\circle*{1.50}}
\put(34.00,14.00){\circle*{1.50}}
\put(9.00,9.00){$1$}
\put(17.00,9.00){$2$}
\put(25.00,9.00){$3$}
\put(33.00,9.00){$4$}
\put(45.00,14.00){$\rightarrow$}
\put(55.00,14.00){$R(\vec{m},\vec{\tau})\sim$}
\put(85.00,20.00){\line(1,0){16.00}}
\put(85.00,20.00){\line(0,-1){6.00}}
\put(101.00,20.00){\line(0,-1){6.00}}
\put(85.00,14.00){\circle*{1.50}}
\put(93.00,14.00){\circle*{1.50}}
\put(101.00,14.00){\circle*{1.50}}
\put(109.00,14.00){\circle*{1.50}}
\put(84.00,9.00){$1^{1,1}$}
\put(92.00,9.00){$2^{1,2}$}
\put(100.00,9.00){$3^{1,1}$}
\put(108.00,9.00){$4^{1,1}$}
\end{picture}
$\;$\\
The corresponding moments are given by
$$
(i) \;\;\phi_{1}(x_{1}x_{3})\phi_{2}(x_{2}x_{4}),\;\;
(ii)\;\;\phi_{1}(x_{1}x_{3})\phi_{2}(x_{2})\phi_{2}(x_{4}),
$$
respectively. We can see how the filters make certain connections
in the partition $R$ disappear.\\[10pt]
\myownsection
\begin{center}
{\sc 5. Convolution limit theorems}
\end{center}
In this section we will prove the central limit theorem
and Poisson's limit theorem for filtered convolutions of states on 
the bialgebra $\widehat{\cal B}$. We choose the convolution formulation
for clarity of exposition, but the general case, based on 
the filtered product of states, is done in an analogous fashion. 

The combinatorics of filtered convolution powers 
$$
\widehat{\phi}^{\star N}=
\widehat{\phi}^{\otimes N}\circ \widehat{\Delta}^{N-1}
$$
where $N\in {\bf N}$, is given by Lemma 5.1. 
To a large extent we follow our approach for the convolution powers
of $q$-deformed states on $U_{q}(su(2))$ given in [10].\\
\indent{\par}
{\sc Lemma 5.1.}
{\it Let $\vec{k}=(k_{1}, \ldots , k_{n})$,
$\vec{\sigma}=(\sigma_{1}, \ldots , \sigma_{n})$, where $k_{i}\in {\bf N}$,
$\sigma_{i}\in {\cal P}({\bf N})$, $1\leq i \leq n$, and let
$N\in {\bf N}$. Then}
\begin{equation}
\label{5.1}
{\widehat\phi}^{\star N}(X_{k_{1}}(\sigma_{1})\ldots X_{k_{n}}(\sigma_{n}))
=\sum_{p=1}^{n}(N)_{p}
\sum_{R=\{R_{1}, \ldots , R_{p}\}\in {\cal P}_{n}}
\prod_{B \in R(\vec{k},\vec{\sigma})}\widehat{\phi} (X_{B})
\end{equation}
{\it where $(N)_{p}=N(N-1)\ldots (N-p+1)$ and 
$X_{B}=\prod_{i\in B}X_{k_{i}}(\sigma_{i})$ for the block $B$ 
of the partition $R(\vec{k}, \vec{\sigma})$, with the product taken in the 
natural order.}\\[5pt]
{\it Proof.}
Denote $X_{1}=X_{k_{1}}(\sigma_{1}), \ldots , X_{n}=X_{k_{n}}(\sigma_{n})$.
Using the notation of (2.11), we have
$$
\widehat{\phi}^{\star N}(X_{1}\ldots X_{n})
=
\sum_{l_{1},\ldots , l_{n} =1}^{N}
\widehat{\phi}^{\otimes N}
(\widehat{j}_{l_{1},N}(X_{1})\ldots \widehat{j}_{l_{n},N}(X_{n}))
$$
and
$$
\widehat{j}_{l_{1},N}(X_{1})\ldots \widehat{j}_{l_{n},N}(X_{n})=
\prod_{m=1}^{n}P(\sigma_{m})^{\otimes (l_{m}-1)}\otimes X_{m}\otimes 
P(\sigma_{m})^{\otimes (N-l_m)}.
$$
The tuple $(l_{1}, \ldots , l_{n})$ defines a partition
$R$ of the set $\{1, \ldots , n\}$ in the usual way. Namely, if
$\{l_{1}, \ldots ,l_{n}\}=\{k_{1}, \ldots , k_{r}\}$, where
$k_{j}$'s are all different, then $R_{j}=\{i: l_{i}=k_{j}\}$.
Thus, from (2.15) we get
$$
\widehat{\phi}^{\otimes N}
(\widehat{j}_{l_{1},N}(X_{1})\ldots \widehat{j}_{l_{n},N}(X_{n}))
=
\prod_{i=1}^{r}\widehat{\phi}(\xi_{i}^{R}(X_{1}\ldots X_{n}))
$$
where $\xi_{i}^{R}$ is a multiplicative extension of the mapping
$$
\xi_{i}^{R}(X_{p})=
\left\{
\begin{array}{ccc}
P(\sigma_{p}) & {\rm if} & p\notin R_{i}\\
X_{p} & {\rm if} & p\in R_{i}
\end{array}
\right..
$$
Now, if $1\leq r \leq n$, then for each partition $R$ consisting of
$r$ blocks, there are $(N)_{r}$ tuples 
$(l_{1}, \ldots , l_{n})$ which give the same contribution 
$\prod_{i}^{r}\widehat{\phi}(\xi_{i}^{R}(X_{1}\ldots X_{n}))$
(the same combinatorial argument is presented in [10] in more detail).
Thus
$$
\widehat{\phi}^{\star N}(X_{1}\ldots X_{n})
=\sum_{r=1}^{n}(N)_{p}
\sum_{R=\{R_{1}, \ldots , R_{p}\}}\prod_{i=1}^{r}\widehat{\phi}
(\xi_{i=1}^{R}(X_{1}\ldots X_{n})).
$$
Finally, note that
$$
\prod_{j=1}^{r}\widehat{\phi}(\xi_{i}^{R}(X_{1}\ldots X_{n}))=
\widehat{\phi}(X_{B_{1}})\ldots \widehat{\phi}(X_{B_{r}})
$$
where $B_{1},\ldots , B_{r}$ are blocks of the partition 
$R(\vec{k}, \vec{\sigma})$ since every block $R_{j}$ of $R$
splits up into subblocks for which all $k_{i}$'s are the same
and are not separated by any filters due to the way $\widetilde{\phi}$
separates words. It is also worth noting that 
$\widehat{\phi}(X_{B})=\phi(X^{\#B})$ where $\#$ stands for the number
of elements.
\hfill $\Box$\\
\indent{\par}
In order to state the central limit theorem, let us introduce the gradation
on $\widehat{\cal B}$ given by $d(X_{k}(\sigma))=1$ and $d(P(\sigma))=0$
for all $k$ and $\sigma$. Then, for $N\in {\bf N}$, define 
$$
D_{1/\sqrt{N}}(W)=\frac{1}{N^{d(W)/2}}W
$$
where $W$ is a word in $\widehat{\cal B}$ and $d(W)$ is its
degree. \\
\indent{\par}
{\sc Corollary 5.2.} 
{\it Consider a family of states $\phi_{N}$ on ${\bf C}[Y]$, where
$N\in {\bf N}$ and suppose that the limits}
$$
\lim_{N\rightarrow \infty}\phi_{N}(Y^{k})=Q(k)
$$
{\it exist and are finite for all $k\in {\bf N}$. Then}
\begin{equation}
\label{5.2}
\lim_{N\rightarrow \infty}\widehat{\phi}^{\star N}_{N}
(X_{k_{1}}(\sigma_{1})\ldots X_{k_{n}}(\sigma_{n}))
=
\sum_{R\in {\cal P}_{n}(\vec{k},\vec{\sigma})}\prod_{B\in R}
Q(\# B)
\end{equation}
{\it where $\# B$ is the number of elements in the block $B$.}\\[5pt] 
{\it Proof.}
It is an immediate consequence of Lemma 5.1 since if 
$R\in {\cal P}_{n}\setminus {\cal P}_{n}(\vec{k}, \vec{\sigma})$, then
the number of blocks in $R(\vec{k}, \vec{\sigma})$
is strictly greater than the number of blocks in $R$
which makes the contribution from $R$ disappear
as $N\rightarrow \infty$. \hfill $\Box$\\
\indent{\par}
{\sc Theorem 5.3.} ({\sc Central limit theorem})
{\it Let $k_{i}\in {\bf N}$, $\sigma_{i}\in {\cal P}({\bf N})$, 
$i=1, \ldots ,n$.
Suppose that $\widehat{\phi}(X_{k_{i}}(\sigma_{i}))=0$ and 
$\widehat{\phi}(X_{k_{i}}^{2}(\sigma_{i}))=1$ for 
$i=1, \ldots , n$. If $n$ is even, then}
\begin{equation}
\label{5.3}
\lim_{N\rightarrow \infty}
\widehat{\phi}^{\star N}\circ D_{1\sqrt{N}}
(X_{k_{1}}(\sigma_{1})\ldots X_{k_{n}}(\sigma_{n}))
=
|{\cal P}_{n}^{\rm pair}(\vec{k}, \vec{\sigma})|
\end{equation}
{\it and, if $n$ is odd, the limit vanishes.}\\[5pt]
{\it Proof.}
It is enough to use Lemma 5.1 and notice that if there is a singleton
in $R$, then there is no contribution from such a partition to 
the right hand side of (5.1). In turn, if there are no singletons, then 
$(N)_{p}/N^{n/2}\rightarrow 0$ unless $2p=n$. That means that
in the limit only pair-partitions may give a non-zero contribution. 
However, note that those pair partitions 
which are not $(\vec{k}, \vec{\sigma})$-adapted
give zero since in that case the number of blocks of 
$R(\vec{k},\vec{\sigma})$ is strictly greater than the number of 
blocks of $R$ and $\prod_{B\in R(\vec{k}, \vec{\sigma})}\widehat{\phi}
(X_{B}) =0$ by the mean zero assumption.\hfill $\Box$ \\
\indent{\par}
{\it Example 1.}
Note that if $\sigma_{i}={\bf N}$ for all $1\leq i \leq n$ and all $n$,
we obtain $|{\cal P}_{n}^{\rm pair}(\vec{k})|$ on the RHS of (5.3) which gives
the moments of the classical multivariate Gaussian law.\\
\indent{\par}
{\it Example 2.}
Here we give some one-dimensional examples. 
If $k_{i}=k$ and $\sigma_{i}=\sigma$ for $1\leq i \leq n$ and all $n$, 
then we obtain
the Gaussian law if $k\in \sigma$ and the $1$-free (or, Boolean)
central limit law corresponding to the discrete measure 
$\mu^{(1)}=1/2(\delta_{-1} +\delta_{1})$ if $k\notin \sigma$.
In turn, if we take 
\begin{equation}
\label{5.4}
\Delta^{(m)}=\widehat{\Delta} \circ \bar{i}^{(m)}
\end{equation}
where $\bar{j}^{(m)}$ is given by (3.3), 
we obtain the $m$-free coproduct defined in [11], for which
the convolution powers tend to the $m$-free central limit laws
and approximate pointwise the Wigner semi-circle law for 
$m=\infty$. For details, see [4].\\
\indent{\par}
{\sc Theorem 5.4.} ({\sc Poisson's limit theorem})
{\it Under the assumptions of Corollary 5.2, 
suppose that $Q(k)=\lambda$ for all $k\in {\bf N}$, where $\lambda >0$. 
Then}
\begin{equation}
\label{5.5}
\lim_{N\rightarrow \infty}
\widehat{\phi}^{\star N}_{N}(X_{k_{1}}(\sigma_{1}) \ldots 
X_{k_{n}}(\sigma_{n}))=
\sum_{R\in {\cal P}_{n}(\vec{k},\vec{\sigma})}
\lambda^{b(R)}
\end{equation}
{\it where $b(R)$ is the number of blocks of $R$.}\\[5pt]
{\it Proof.}
It is an immediate consequence of Corollary 5.2.\hfill $\Box$\\
\indent{\par}
{\it Example 1.}
Let us first give some one-dimensional examples.
Again, if $k_{i}=k$ and $\sigma_{i}=\sigma$ for $1\leq i \leq n <\infty$, then
we obtain the classical Poisson law for $k\in \sigma$ and the 1-free
(or Boolean) Poisson law for $k\notin \sigma$ corresponding to
the discrete measure $\mu_{\lambda}^{(1)}=1/(1+\lambda)
(\delta_{0}+\lambda \delta_{1+\lambda})$. Considering
linear combinations of sample sums as in the preceeding
example, we obtain the $m$-free Poisson laws 
for $m\in {\bf N}$ 
and the free Poisson law [21] for $m=\infty$ (see [4]).\\
\indent{\par}
{\it Example 2.}
If we take $\widehat{\phi}$ given by (2.16), i.e. corresponding to the
product of measures, then we can generalize Lemma 5.1 and Corollary 5.2
to the effect that if $\lim_{N\rightarrow\infty}\phi_{s,N}(Y^{k})=\lambda_{s}$
where $\lambda_{s}>0$, $s\in {\bf N}$, then the RHS of (5.5) takes the form
$$
\sum_{R\in {\cal P}_{n}(\vec{k}, \vec{\sigma})}\lambda_{s_{1}}
\lambda_{s_{2}}\ldots \lambda_{s_{p}}
$$
where $s_{1}, \ldots , s_{p}$ correspond to the blocks $B_{1}, \ldots ,B_{p}$
of the partition $R(\vec{k}, \vec{\sigma})$ and denote their colors
(which are the same within one block by (A1) of Definition 4.1). These moments
are the moments of the multivariate classical Poisson law.\\[10pt]
\myownsection
\begin{center}
{\sc 6. Filtered Fundamental Operators}
\end{center}
In this section we recall basic facts concerning multiple 
symmetric Fock spaces over ${\cal K}\equiv L^{2}({\bf R}^{+})$, which will be
the underlying space for the filtered fundamental processes.

Let ${\cal G}$ be a separable Hilbert space with a countably infinite 
fixed orthonormal basis
$(e_{n})_{n\in {\bf N}}$. It is called the {\it multiplicity space}.
By a {\it multiple symmetric Fock space } over
${\cal K}$ we understand the symmetric Fock space over
${\cal H}=L^{2}({\bf R}^{+}, 
{\cal G})\cong L^{2}({\bf R}^{+})
\otimes {\cal G}\equiv {\cal K}\otimes {\cal G}$, namely
$$
\Gamma({\cal H})={\bf C}\Omega \oplus \bigoplus_{n=1}^{\infty}{\cal H}^{\circ n}
$$
where ${\cal H}^{\circ n}$ denotes the $n$-th symmetric tensor power
of ${\cal H}$ and $\Omega$ is the vacuum vector, 
with the scalar product given by $\langle\Omega , \Omega \rangle=1$,
$\langle \Omega , u \rangle =0$ and
$$
\langle u_{1} \circ \ldots \circ u_{n},
v_{1} \circ \ldots \circ v_{m} \rangle
=\delta_{n,m}\frac{1}{n!}
\sum_{\sigma \in {\bf S}_{n}}
\langle u_{1} , v_{\sigma(1)} \rangle
\ldots 
\langle u_{n} , v_{\sigma(n)} \rangle
$$
where
$$
u_{1} \circ \ldots \circ u_{n}=
\frac{1}{n!}\sum_{\sigma \in {\bf S}_{n}}
u_{\sigma(1)}\otimes \ldots \otimes u_{\sigma(n)}
$$
and ${\bf S}_{n}$ denotes the symmetric group of order $n$.

Denote by ${\cal H}^{(\sigma)}$ the linear subspace of ${\cal H}$ spanned by
all $u\in {\cal H}$ of the form
$$
u=\sum_{k\in \sigma}u^{(k)}\otimes e_{k}
$$
where $\sigma\in {\cal P}({\bf N})$. In particular, 
we put ${\cal H}^{(\emptyset)}=\{0\}$. The set $\sigma$ will be called
a {\it filter} and the associated canonical projection will be denoted
$\Pi^{(\sigma)}\;: {\cal H}\rightarrow {\cal H}^{(\sigma)}$
with $v^{(\sigma)}=\Pi^{(\sigma)}v$
for any $v\in {\cal H}$. Then let
$P^{(\sigma)} \;: \Gamma({\cal H})\rightarrow \Gamma({\cal H}^{(\sigma)})$
be the second quantization of $\Pi^{(\sigma)}$. Thus, if  
$\varepsilon (v)$ is an exponential vector in $\Gamma({\cal H})$,
i.e.
$$
\varepsilon (v)=\bigoplus_{n=0}^{\infty}\frac{1}{\sqrt{n!}}v^{\otimes n}
$$
with $v^{\otimes 0}=\Omega$, we have 
$P^{(\sigma)}\varepsilon(u)= \varepsilon(u^{(\sigma)})$.

Of special importance will be
subspaces ${\cal H}^{(r)}$ of ${\cal H}$ spanned by
all $u\in {\cal H}$ of the form
$$
u=\sum_{k=1}^{r-1}u^{(k)}\otimes e_{k},
$$ 
where $r >1$, i.e. here $\sigma=\{1, \ldots , r-1\}$; we set 
${\cal H}^{(1)}=\{0\}$.
In $\Gamma({\cal H})$, we will use the {\it finite particle domain} 
$\Gamma_{0}({\cal H})$, i.e. the linear space generated by vectors
of the form
$$
v_{1}\circ v_{2} \circ \ldots \circ v_{n}
$$
where $v_{1}, \ldots , v_{n}\; \in {\cal H}$, $n\in {\bf N}$.

Since ${\cal H}$ can be viewed as a direct sum of infinitely
many copies of ${\cal K}$ and we need some convenient
terminology concerning the numbering of those copies, we will
refer to them as {\it colors}. Thus, in the direct sum
decomposition 
$$ 
{\cal H}=\bigoplus_{k\in {\bf N}} {\cal K}\otimes e_{k} 
$$
the $k$-th summand will be associated with 
the $k$-th color and we will say that 
non-zero vectors from that summand are of $k$-th color. In addition,
to the zero vector we assign the $0$-th color.

By {\it filtered creation and 
annihilation operators} we will understand operators given by
\begin{eqnarray}
\label{6.1}
a^{(\sigma)*}(f\otimes e_{k})&=&a^{*}(f\otimes e_{k})P^{(\sigma)}\\
\label{6.2}
a^{(\sigma)}(f\otimes e_{k})&=&P^{(\sigma)}a(f\otimes e_{k}),
\end{eqnarray}
respectively, where $a^{*}(f\otimes e_{k})$ and $a(f\otimes e_{k})$
are the usual boson creation and annihilation operators (see [P]).
Thus, filtered
creation operators first ``filter out particles of colors 
which are not in $\sigma$ and then create a particle of given color'', 
whereas the filtered annihilation operators 
``first annihilate a particle of a given color and then filter
out particles of colors which are not in $\sigma$''.

In addition, we define 
\begin{equation}
\label{6.3}
a^{(k,\sigma)\circ}
=
a^{(k)\circ}P^{(\sigma\cup \{k\})} 
\end{equation}
and call {\it filtered number operators}.
In an analogous fashion one can define exchange operators. \\
\indent{\par}
{\sc Proposition 6.1.}
{\it The finite particle domain $\Gamma_{0}({\cal H})$ 
is contained in the domains of filtered fundamental operators.
Furthermore, the following relations hold:}
\begin{eqnarray*}
a^{(\sigma)}(f\otimes e_{k})
(v_{1}\circ v_{2} \circ \ldots \circ v_{n})
&=&\frac{1}{\sqrt{n}}
\sum_{j=1}^{n}
\langle f , v^{(k)}_{j} \rangle 
v_{1}^{(\sigma)} \circ \ldots \circ \breve{v}_{j} \circ
\ldots \circ v_{n}^{(\sigma)},\\[5pt]
a^{(\sigma)*}(f\otimes e_{k})
(v_{1}\circ v_{2} \circ \ldots \circ v_{n})
&=&
\sqrt{n+1}(f\otimes e_{k})\circ v_{1}^{(\sigma)}\circ \ldots \circ v_{n}^{(\sigma)},\\[5pt]
a^{(k,\sigma)\circ} (v_{1}\circ v_{2}\circ \ldots \circ v_{n})
&=&
\sum_{j=1}^{n}v_{1}^{(\sigma)}\circ \ldots \circ (v_{j}^{(k)}\otimes e_{k})
\circ \ldots \circ v_{n}^{(\sigma)}
\end{eqnarray*}
{\it with 
$a^{(\sigma)}(f\otimes e_{k})\Omega =0$, 
$a^{(\sigma)*}(f\otimes e_{k})\Omega = f\otimes e_{k}$ and
$a^{(k,\sigma)\circ} \Omega = 0$,
where $v_{1}, \ldots , v_{n}\in {\cal H}$, 
$k\in {\bf N}$, $\sigma\in {\cal P}({\bf N})$,
$n\in {\bf N}$.}\\[5pt]
{\it Proof.}
The first statement follows from the definitions (6.1)-(6.3) and
an analogous property of the canonical
(CCR) operators and the fact that the projections $P^{(\sigma)}$ leave the finite
particle domain invariant.
Similarly, the relations follow immediately from the analogous formulas 
for the canonical (CCR) operators (we use Hudson-Parthasarathy's normalization).
\hfill
$\Box$\\
\indent{\par}
{\sc Lemma 6.2.}
{\it Filtered creation and annihilation operators satisfy the following
relations on the finite particle domain:}
$$
a^{(\sigma)}(f\otimes e_{k})a^{(\tau)*}(g\otimes e_{l})-
a^{(\tau)*}(g\otimes e_{l})a^{(\sigma)}(f\otimes e_{k})P^{(\tau)}\1_{\{l\in \sigma\}}
=\delta_{k,l}\langle f , g \rangle P^{(\sigma \cap \tau)}
$$
{\it for any $k,l \in {\bf N}$, $\sigma,\tau \in {\cal P}({\bf N})$,
$f,g \in {\cal K}$}.\\[5pt]
{\it Proof.}
In the proof given below we understand that the
equations hold on the finite particle domain, but
it remains valid on the whole intersection of the domains 
of the considered filtered operators.
Using canonical commutation relations (CCR) of the form
$$
a(f\otimes e_{k})a^{*}(g\otimes e_{l})-a^{*}(g\otimes e_{l})a(f\otimes e_{k})
=\delta_{k,l}\langle f , g \rangle.
$$
we obtain
$$
a^{(\sigma)}(f\otimes e_{k})a^{(\tau)*}(g\otimes e_{l})=
P^{(\sigma)}a(f\otimes e_{k})a^{*}(g\otimes e_{l})P^{(\tau)}
$$
$$
=P^{(\sigma)}a^{*}(g\otimes e_{l})a(f\otimes e_{k})P^{(\tau)}+
\delta_{k,l}\langle f , g \rangle P^{(\sigma\cap \tau)}.
$$
Now, note that if $l\notin \sigma$, then
$$
P^{(\sigma)}a^{*}(g\otimes e_{l})a(f\otimes e_{k}) P^{(\tau)}=0
$$
and we obtain
$$
a^{(\sigma)}(f\otimes e_{k})a^{(\tau)*}(g\otimes e_{l})=\delta_{k,l}\langle f ,
g \rangle P^{(\sigma \cap \tau)}.
$$
Consider now the case $l\in \sigma$. Then
$$
P^{(\sigma)}a^{*}(g\otimes e_{l}) a(f\otimes e_{k}) P^{(\tau)}-
a^{*}(g\otimes e_{l})P^{(\sigma \cap \tau)}a(f\otimes e_{k})P^{(\tau)}
$$
$$
=a^{*}(g\otimes e_{l})(P^{(\sigma)}-P^{(\sigma \cap \tau)})a(f\otimes e_{k})P^{(\tau)}
$$
since $P^{(\sigma)}$ commutes with $a^{*}(g\otimes e_{l})$ for $l\in \sigma$.
Note that if $\sigma \subseteq \tau$, then $P^{(\sigma)}=P^{(\sigma \cap \tau)}$ 
and thus the above expression vanishes. In turn, if $\tau \subseteq \sigma$, then 
$P^{(\sigma)}-P^{(\sigma \cap \tau)}= P^{(\sigma)}-P^{(\tau)}$. However,
$$
P^{(\tau)}:\;\; \Gamma({\cal H})\rightarrow \Gamma({\cal H}^{(\tau)})
$$ 
and $a(f\otimes e_{k})$ leaves $\Gamma({\cal H}^{(\tau)})$ invariant, hence
when we apply $P^{(\sigma)}-P^{(\tau)}$, we can see that the above
expression also vanishes. Therefore, if $l\in \sigma$, we obtain
$$
a^{(\sigma)}(f\otimes e_{k})a^{(\tau)*}(g\otimes e_{l})-
a^{*}(g\otimes e_{l})P^{(\sigma \cap \tau)}a(f\otimes e_{k})P^{(\tau)}
$$
$$
=
a^{(\sigma)}(f\otimes e_{k})a^{(\tau)*}(g\otimes e_{l})-
a^{(\tau)*}(g\otimes e_{l})a^{(\sigma)}(f\otimes e_{k})P^{(\tau)}
$$
$$
=
\langle f , g \rangle \delta_{k,l} P^{(\sigma \cap \tau)}.
$$
Combining the two cases $l\in \sigma$ and $l \notin \sigma$  ends the proof.
\hfill $\Box$\\
\indent{\par}
Let us finally define the fundamental procesess associated
with the filtered fundamental operators introduced in this section.
They will appear in Sections 5-7 when finding GNS realizations
of limit states. They will also serve as integrators in the filtered calculus
developed in [8].
Thus, in connection with (6.1)-(6.3), let
\begin{eqnarray}
\label{6.4}
A^{(k,\sigma)*}_{t}&=&a^{(\sigma)*}(\chi_{[0,t]}\otimes e_{k}),\\
\label{6.5}
A^{(k,\sigma)}_{t}&=&a^{(\sigma)}(\chi_{[0,t]}\otimes e_{k}),\\
\label{6.6}
A^{(k,\sigma)\circ}_{t}&=&
\lambda(I_{[0,t]}\otimes |e_{k}\rangle \langle e_{k}|)
P^{(\sigma\cup \{k\})},\\
\label{6.7}
A^{(0,\sigma)}_{t}&=& tP^{(\sigma)}
\end{eqnarray}
where $t\geq 0$, $k\in {\bf N}$, $\sigma\in {\cal P}({\bf N})$,
$I_{[0,t]}$ denotes the
operator of mulitplication by the characteristic function 
$\chi_{[0,t]}$ on $L^{2}({\bf R}^{+})$, and $\lambda (H)$
denotes the differential second quantization of $H\in {\cal B}({\cal H})$. 
The 
families of processes given by (6.4)-(6.7)
will be called {\it filtered creation, annihilation, 
number} and {\it time procesess}, respectively.
When speaking of all of them, we 
will call them {\it filtered fundamental processes}.
By {\it filtered Brownian motion} we will understand
the unital *-algebra generated by filtered 
creation and annihilation operators indexed by time intervals.\\[10pt]
\myownsection
\begin{center}
{\sc 7. Random walk on the filtered bialgebra}
\end{center}
In this section we show that a limit of continuous-time random walks on 
the filtered *-bialgebra gives the filtered Brownian motion.
This gives a multivariate Brownian motion on the multiple symmetric Fock space
which satisfies the properties required by the axioms for white 
noise on *-bialgebras given in [1] and [18] and includes quantum Brownian
motions for different types of independence [2]. For the first quantum
version of the Wiener process, see [3]. We follow the notation used in [9]
for the random walk on $U_{q}(su(2))$.

Instead of $\widehat{\cal B}$, we choose to work with a slightly more general 
unital *-bialgebra $\widehat{\cal C}$, also called filtered *-bialgebra,
which is defined to be the unital *-algebra over ${\bf C}$ generated
by $X_{k}(\sigma)$, $X_{k}^{*}(\sigma)$ and $P(\sigma)$,
where $k\in {\bf N}$, 
$\sigma\in {\cal P}({\bf N})$ subject to relations (2.6)-(2.7), where
$P(\sigma)$ is a projection for each $\sigma\in {\cal P}({\bf N})$
(this of course also means that $P(\sigma) $ commutes with $X_{k}^{*}(\tau)$ 
for $k\in \sigma$), with the coproduct in which
$X_{k}(\sigma)$ and $X_{k}^{*}(\sigma)$ are both $P(\sigma)$-primitive
and $P(\sigma)$ is group-like (cf. (2.8)-(2.9)).

Let $k\in {\bf N}$, $\sigma \in {\cal P}({\bf N})$, $N\in {\bf N}$, and
consider a sequence of continuous-time random walks on $\widehat{\cal C}$ 
given by 
\begin{eqnarray}
\label{7.1}
\widehat{\Delta}_{s,t}^{N}(X_{k}^{\natural}(\sigma))
&=&
\sum_{l=N_{s}+1}^{N_{t}}
P(\sigma)^{\otimes (l-1)}\otimes 
X_{k}^{\natural}(\sigma) \otimes P(\sigma)^{\otimes \infty}\\
\label{7.2}
\widehat{\Delta}_{s,t}^{N}(P(\sigma))&=&
1^{\otimes N_{s}}\otimes P(\sigma)^{\otimes (N_{t}-N_{s})}\otimes 
1^{\otimes \infty}
\end{eqnarray}
where $0\leq s \leq t <\infty$, $N_{t}=E[tN]$, 
with $X_{k}^{\natural}(\sigma)\in
\{X_{k}(\sigma), X_{k}^{\natural}(\sigma)\}$.

It is easy to check 
that for each pair $(s,t)$ and natural number $N$, the mapping
$$
\widehat{\Delta}_{s,t}^{N}: 
\;\widehat{\cal C}\rightarrow \widehat{\cal C}^{\otimes \infty}
$$ 
given by the linear and multiplicative extension of (7.1)-(7.2)
is a unital *-homomorphism and the triple
$(\widehat{\cal C}^{\otimes \infty}, (\Delta_{s,t}^{N})_{0\leq s\leq t},
\widehat{\phi}^{\otimes \infty})$
satisfies for each $N\in {\bf N}$
the properties required from a stochastic process over
the bialgebra $\widehat{\cal C}$ given in [1]. In particular,
$$
\Delta_{s,t}^{N}\star \Delta^{N}_{t,r}=\Delta^{N}_{s,r}
$$
for all $0\leq s \leq t \leq r$, where $\Delta_{s,t}^{N}
\star \Delta_{t,r}^{N}=M\circ (\Delta_{s,t}^{N}\otimes \Delta_{t,r}^{N})
\circ \Delta$ with $M(a\otimes b)=ab$. For details on stochastic
processes over *-bialgebras see [1] and [18].

In this *-bialgebra formulation,
further preparations are similar to those which lead to 
the central limit theorem. Namely, for a given state $\phi$ on 
${\bf C}\langle Y, Y^{*} \rangle$ we denote by $\widetilde{\phi}$
its Boolean extension to ${\bf C}\langle Y, Y^{*} , P\rangle$, 
where $P$ is a projection and we set 
$\widehat{\phi}=\widetilde{\phi}^{\otimes \infty}\circ \eta$,
where
$$
\eta :\; \widehat{\cal C}\rightarrow \bigotimes_{k=1}^{\infty}
{\bf C}\langle Y_{k},Y_{k}^{*}, P_{k}\rangle
$$
is defined by the *-multiplicative extension of formulas (2.13)-(2.14).

Below we will study the limit of distributions of the mixed
moments of (7.1) as $N\rightarrow \infty$ and 
find the GNS representation of the limit state.
Let us remark that 
more general sample sums indexed by $f\in L^{2}_{c}({\bf R}^{+})$
can also be given and the proofs of this section will still hold.
\\
\indent{\par}
{\sc Theorem 7.1.}
{\it Let $Z_{i}=X_{k_{i}}(\sigma_{i})$,
$Z_{i}^{*}=X_{k_{i}}^{*}(\sigma_{i})$, where
$k_{i}\in {\bf N}$, $\sigma_{i}\in {\cal P}({\bf N})$,
$i=1, \ldots ,n$. 
Suppose that $\widehat{\Phi}=\widehat{\phi}^{\otimes \infty}$ with
$\phi(Y)=0$ and the only non-vanishing second-order 
moment of $\phi$ is given by $\phi(YY^{*})=1$. Then}
$$
\lim_{N\rightarrow \infty}
N^{-n/2}
\widehat{\Phi} (\widehat{\Delta}_{s_{1},t_{1}}^{N}
(Z_{1}^{\natural})\ldots \widehat{\Delta}_{s_{n},t_{n}}^{N}(Z_{n}^{\natural}))
=
\varphi (a^{(\sigma_{1})\natural}(v_{1})\ldots
a^{(\sigma_{n})\natural}(v_{n}))
$$
{\it where $Z_{i}^{\natural} \in\{Z_{i}, Z_{i}^{*}\}$,
$v_{i}=\chi_{[s_{i},t_{i}]}\otimes e_{k_{i}}$,
$i=1, \ldots , n$, and
$\varphi (.)=\langle \Omega , . \Omega \rangle$ is
the vacuum expectation in $\Gamma({\cal H})$.}\\[5pt]
{\it Proof.}
From the general invariance principle [22] and the combinatorics
of the filtered central limit theorem it follows that
for even $n=2p$ we have
$$
LHS
=
\sum_{R\in {\cal P}_{2p}^{{\rm pair}}(\vec{k},\vec{\sigma})}
\delta(R)
\prod_{m=1}^{p}
\langle 
v_{\alpha (m)},v_{\beta (m)}
\rangle
$$
where the blocks of $R$ consist of two-element sets
$\{\alpha(i), \beta(i)\}$, 
with $\alpha(i)<\beta(i)$, $i=1, \ldots , p$ and $\delta(R)=1$ if
for the given partition $R$ we have
$Z_{\alpha(i)}^{\natural}=Z_{\alpha(i)}$ and $Z_{\beta(i)}^{\natural}=
Z_{\beta(i)}^{*}$ and otherwise $\delta(R)=0$.
It is clear that if $n$ is odd, then $LHS=0$.

It is clear that $RHS=0$ if $n$ is odd, too -- 
it is enough to use the properties of creation and 
annihilation operators following from Proposition 6.1.
Therefore assume that $n=2p$.
Next, notice that in order that $LHS=RHS$ it is enough to show
the following claim:
$$
\varphi 
(a_{1}^{\natural}\ldots
a_{n}^{\natural})=
\sum_{R\in 
{\cal P}_{2p}^{{\rm pair}}(\vec{k},\vec{\sigma})}
\delta(R)
\prod_{m=1}^{p}
\langle 
v_{\alpha (m)}, v_{\beta (m)}
\rangle  
$$
where, for simplicity, we denote 
$a_{j}^{\natural}=a^{\sigma_{j}\natural}(v_{j})$. 

We begin with the simplest case, i.e. 
\begin{eqnarray*}
\varphi
(a_{1}a_{2}^{*}\ldots a_{2p-1}a_{2p}^{*})&=&
\langle v_{1}, v_{2} \rangle \langle v_{3}, v_{4} \rangle \ldots
\langle v_{2p-1} , v_{2p} \rangle \\
&=&
\sum_{R\in {\cal P}_{2p}^{{\rm pair}}(\vec{k},\vec{\sigma})}
\delta(R)
\prod_{m=1}^{p}
\langle 
v_{\alpha (m)}, v_{\beta (m)}
\rangle ,
\end{eqnarray*}
the second expression being formally written as a sum since we have at most
one partition contributing to it. 
This is the beginning of an induction procedure. 
Namely, it is enough to show that from the claim 
being true for all expectations
of orders $\leq 2p-2$ and for
$\varphi(a_{1}^{\natural}
\ldots a_{i}^{*}a_{i+1}\ldots a_{2p}^{\natural})$
it follows that it also holds for
the expectation of the form 
$\varphi(a_{1}^{\natural}\ldots 
a_{i+1}a_{i}^{*}\ldots a_{2p}^{\natural})$.

Suppose that $k_{i}\notin \sigma_{i+1}$. Then
\begin{eqnarray*}
\varphi(a_{1}^{\natural}\ldots 
a_{i+1}a_{i}^{*}\ldots a_{2p}^{\natural})
&=&\langle v_{i}, v_{i+1} \rangle
\varphi
(a_{1}^{\natural}\ldots 
P^{(\sigma_{i}\cap \sigma_{i+1})} \ldots
a_{2p}^{\natural})\\
&=&
\sum_{R\in {\cal P}_{2p}^{{\rm pair}}(\widehat{k},\widehat{\sigma}|i,i+1)}
\delta(R)
\prod_{m=1}^{p}
\langle 
v_{\alpha (m)}, v_{\beta (m)}
\rangle  
\end{eqnarray*}
where ${\cal P}_{n}^{{\rm pair}}(\widehat{k},\widehat{\sigma}|i,i+1)$ 
denotes all $(\widehat{k}, \widehat{\sigma})$-- adapted pair 
partitions with
$$
\widehat{k}=(k_{1}, \ldots , k_{i+1}, k_{i}, \ldots ,k_{n}),\;\;
\widehat{\sigma}=(\sigma_{1}, \ldots , \sigma_{i+1}, \sigma_{i}, \ldots , \sigma_{n})
$$
in which $(i,i+1)$ forms a pairing. The first equality follows
from filtered relations of Lemma 6.2, 
whereas the second -- from the inductive assumption
and the fact that if $k_{i}\notin \sigma_{i+1}$, then
$$
\forall R\in {\cal P}_{2q}^{\rm pair}(\widehat{k},\widehat{\sigma})\;\; 
\exists R' \in {\cal P}_{2q-2}^{\rm pair}(\zeta_{i,i+1}(\vec{k},\vec{\sigma}))\;:\;\;
R=R'\cup \{(i,i+1)\},
$$
where $\zeta_{i,i+1}(\vec{k},\vec{\sigma})=(\widetilde{k}, \widetilde{\sigma})$, with
$$
\widetilde{k}= (k_{1}, \ldots , k_{i-1},k_{i+2}, \ldots , k_{n}),\;\;
\widetilde{\sigma}= 
(\sigma_{1}, \ldots , \sigma_{i-1},\zeta(\sigma_{i+2}), \ldots , \zeta(\sigma_{n})),
$$
and 
$$
\zeta(\sigma_{l})=\left\{
\begin{array}{cc}
\sigma_{l}\cap \sigma_{i}\cap \sigma_{i+1} 
& {\rm if}\;\;(s,l) \;\;{\rm is}\; {\rm a}\;{\rm pairing}\; {\rm for}
\;\; s<i\\
\sigma_{l} & {\rm otherwise}
\end{array}
\right. 
$$
where $i+2\leq l \leq n$.

In turn, if $k_{i}\in \sigma_{i+1}$, then using Lemma 6.2 again, we obtain
$$
\varphi(a_{1}^{\natural}\ldots a_{i+1}a_{i}^{*}\ldots a_{2p}^{\natural})
$$
$$
=\langle v_{i+1}, v_{i} \rangle
\varphi
(a_{1}^{\natural}\ldots 
P^{(\sigma_{i}\cap \sigma_{i+1})} \ldots
a_{2p}^{\natural})
+
\varphi(a_{1}^{\natural}\ldots 
a_{i}^{*}
a_{i+1}P^{(\sigma_{i})}\ldots 
a_{2p}^{\natural}).
$$
By the inductive assumption and similar arguments as above, 
the first term gives
$$
\sum_{R\in {\cal P}_{2p}^{{\rm pair}}(\widehat{k},\widehat{\sigma}|i,i+1)}
\delta(R)
\prod_{m=1}^{p}
\langle 
v_{\alpha (m)}, v_{\beta (m)}
\rangle  ,
$$
whereas the second,  
a sum over all the remaining partitions 
from ${\cal P}_{n}^{{\rm pair}}(\widehat{k},\widehat{\sigma})$ (it is 
disjoint from the first since $(i,i+1)$
cannot form a pairing as in the associated creation-annihilation
pair the annihilation operator follows the creation operator), namely
$$
\sum_{R\in {\cal P}_{2p}^{{\rm pair}}(\widehat{k},\widehat{\sigma})\setminus 
{\cal P}_{2p}^{{\rm pair}}(\widehat{k},\widehat{\sigma}|i,i+1)}
\delta(R)
\prod_{m=1}^{p}
\langle 
v_{\alpha (m)}, v_{\beta (m)}
\rangle
$$
by the inductive assumption.
Note that the projection $P^{(\sigma_{i})}$, which follows the
annihilation operator in the second term, ensures that the
annihilation operator $a^{(\sigma_{i+1})}(v_{i+1})$
in the original expression cannot be
paired off with any creation operator (standing to the right 
of this annihilation operator) of color $k \notin \sigma_{i}$.

Adding now those two expressions, we can see that the claim
holds for 
$$
\varphi(a_{1}^{\natural}\ldots a_{i+1}a_{i}^{*}\ldots 
a_{2p}^{\natural}),
$$
which finishes the proof.
\hfill $\Box$\\
\indent{\par}
{\it Example.}
If $k_{i}=k$, $\sigma_{i}=\sigma$ for $1\leq i \leq n$ and arbitrary $n$,
then we obtain the CCR Brownian motion if $k\in \sigma $ and 
Boolean Brownian motion
if $k\notin \sigma$. By taking linear combinations of sample sums
corresponding to $m$-free (free) independent random variables, we obtain
the $m$-free (free) Brownian motion [4].\\[10pt]
\myownsection
\begin{center}
{\sc 8. Filtered White Noise}
\end{center}
In this section we define the general notion of 
filtered white noise, determine its
combinatorics and study the example of filtered Poisson
white noises. Our approach largely parallels that used
by Speicher for free white noise [21]. \\
\indent{\par}
{\sc Definition 8.1.}
Let ${\rm Int}({\bf R}^{+})$ denote the intervals in ${\bf R}^{+}$.
An $s$-dimensional {\it filtered white noise} consists of
a unital *-algebra ${\cal C}$, a state $\rho$ on ${\cal C}$ and
a family of finitely additive mappings 
${\rm Int}({\bf R}^{+})\rightarrow {\cal C}$,
$$
I\rightarrow (c_{I}(k,\sigma;1),\ldots , c_{I}(k,\sigma;s)),\;\; k\in {\bf N},\;\;
\sigma\in {\cal P}({\bf N})
$$
such that \\
(i) for any pairwise disjoint intervals $I(1), \ldots , I(n)$,
\begin{equation}
\label{8.1}
\rho (c_{I(l_{1})}(k_{1},\sigma_{1};q_{1})\ldots c_{I(l_{n})}(k_{n},\sigma_{n};q_{n}))
=\rho(c_{B_{1}})\ldots \rho (c_{B_{r}})
\end{equation}
where
$B_{1}, \ldots , B_{r}$ are the blocks of $R(\vec{k},\vec{\sigma})$, with 
$R$ being the partition associated with $(l_{1}, \ldots , l_{n})$ and
$c_{B}$ denotes the product, taken in the natural order, of 
$c_{I(l_{i})}(k_{i},\sigma_{i};q_{i})$'s for $i\in B$,\\[3pt]
(ii)  
the distribution $\rho_{I}=\rho|{\cal C}_{I}$ 
depends only on the Lebesgue measure
of the interval $I$, where ${\cal C}_{I}$ denotes the unital *-algebra
generated by $c_{I}(k,\sigma;q)$, $k\in {\bf N}$, $\sigma\in {\cal P}({\bf N})$ 
and $1\leq q \leq s$.\\
\indent{\par}
{\sc Lemma 8.2}
{\it Let $({\cal C}, \rho , 
(c_{I}(k,\sigma;1), \ldots , c_{I}(k,\sigma;n))
_{I\in {\rm Int}({\bf R}^{+}), k\in {\bf N} ,\sigma\in {\cal P}({\bf N})})$
be an $s$-di\-men\-sio\-nal filtered white noise and let 
$c_{t}(k,\sigma)=c_{[0,t)}(k,\sigma)$ for $k\in {\bf N}$, $\sigma\in {\cal P}({\bf N})$.
Then}
$$
\rho(c_{t}(k_{1},\sigma_{1};q_{1})\ldots c_{t}(k_{n},\sigma_{n};q_{n}))
=\sum_{R\in {\cal P}_{n}(\vec{k},\vec{\sigma})}
\prod_{B\in R}Q_{t}(B)
$$
{\it where}
\begin{eqnarray*}
Q_{t}(B)&=&Q_{t}(k, (\sigma_{i})_{i\in B},(q_{i})_{i\in B})\\
&=&
\lim_{N\rightarrow \infty}\rho 
(
c_{t/N}(k,\sigma_{i(1)};q_{i(1)})
\ldots 
c_{t/N}(k,\sigma_{i(r)};q_{i(s)})
)
\end{eqnarray*}
{\it and $B=\{i(1), \ldots , i(m)\}$ with $i(1)<\ldots <i(m)$
with $k=k(B)=k_{i(j)}$ for all $i(j)\in B$}.\\[5pt] 
{\it Proof.}
By additivity of the filtered white noise, we can split up 
each $c_{t}^{(k,\sigma;q)}$ into a sum of $N$ summands:
$$
c_{t}(k,\sigma;q)=\sum_{l=1}^{N} c_{I(l)}(k,\sigma;q)
$$
where $I(l)=[(l-1)t/N, lt/N)$. Moreover, the summands have the same 
distributions. 

Now, note that from (8.1) it follows that we can use
the same combinatorial argument as in
Corollary 5.2 providing the limits
$$
Q_{t}(B)=\lim_{N\rightarrow \infty}N\rho 
(
c_{I(l)}(k,\sigma_{1};q_{i(1)})\ldots 
c_{I(l)}(k,\sigma_{n};q_{i(r)})
)
$$
$$
=\lim_{N\rightarrow \infty}
N\rho(c_{I(1)}(k,\sigma_{1};q_{i(1)})\ldots c_{I(1)}(k,\sigma_{n};q_{i(r)}))
$$
exist for all $l,k \in {\bf N}$, $i_{1}, \ldots , i_{r}$, $1\leq r \leq n$,
where $k=k(B)$ and the dependence of the limit on 
$k$, $\sigma_{i(1)}, \ldots , \sigma_{i(r)}$ and $q_{i(1)}, \ldots , q_{i(r)}$
is suppressed.
Existence of such limits follows from an induction procedure
which is analogous to that in the free case [21].
\hfill $\Box$\\
\indent{\par}
{\it Example.}
A $2$-dimensional filtered Gaussian noise is obtained from
$c_{t}(k,\sigma;1)=A_{t}^{(k,\sigma)}$, $c_{t}(k,\sigma;2)=A_{t}^{(k,\sigma)*}$,
given by (6.4)-(6.5),
for any $k\in {\bf N}$, $\sigma\in {\cal P}({\bf N})$, $t\geq 0$ with 
$\rho=\varphi$, the vacuum expectation in $\Gamma({\cal H})$. Then
$$
\varphi (c_{t}(k_{1},\sigma_{1};q_{1})\ldots c_{t}(k_{n},\sigma_{n};q_{n}))
=
\left\{
\begin{array}{cll}
\sum_{R\in {\cal P}_{n}^{{\rm pair}}(\vec{k}, \vec{\sigma})}\prod_{B\in R}Q_{t}(B)
& n & {\rm even}\\
0 & n & {\rm odd}
\end{array}
\right.
$$
where the generator $Q_{t}$ does not depend of $k(B)$ and $\sigma_{i}$'s and
is given by
$$
Q_{t}(B)=Q_{t}(i(1),i(2))=
\left\{
\begin{array}{ll}
t& {\rm if}\;\;q_{i(1)}=1,\; q_{i(2)}=2\\
0& {\rm otherwise} 
\end{array}
\right.
$$
Note that the 
filtered Gaussian noise was obtained before as the GNS
representation of the limit state of the invariance principle 
(in a slightly more general version). \\
\indent{\par}
Below we will find the expectations
of the  
filtered (multivariate) Poisson noise constructed from
filtered fundamental processes given by (6.4)-(6.7), 
which is also a way to justify the correctness
of our definition of filtered number operators (6.3). \\
\indent{\par}
{\sc Theorem 8.3.}
{\it For any $k\in {\bf N}$, $\sigma\in {\cal P}({\bf N})$ and $t\geq 0$, let}
\begin{equation}
\label{8.2}
\Lambda_{t}^{(k,\sigma)}=A_{t}^{(k,\sigma)}+A_{t}^{(k,\sigma)*}+A_{t}^{(k,\sigma)}+A_{t}^{(0,\sigma)}
\end{equation}
{\it and let $\varphi$ be the vacuum expectation state in $\Gamma({\cal H})$.
Then}
$$
\varphi(\Lambda_{t}^{(k_{1},\sigma_{1})}\ldots \Lambda_{t}^{(k_{n},\sigma_{n})})
=
\sum_{R\in {\cal P}_{n}(\vec{k}, \vec{\sigma})}
t^{b(R)},
$$
{\it where $k_{1}, \ldots , k_{n}\in {\bf N}$, 
$\sigma_{1}, \ldots , \sigma_{n}\in {\cal P}({\bf N})$ and 
$b(R)$ is the number of blocks of $R$.}\\[5pt]
{\it Proof.}
First of all, notice that if $I(1), \ldots , I(r)$ are disjoint
intervals in ${\bf R}^{+}$, then 
$$
\varphi (\Lambda_{I(l_{1})}^{(k_{1},\sigma_{1})}\ldots 
\Lambda_{I(l_{n})}^{(k_{n},\sigma_{n})})
=
\varphi(\Lambda_{B_{1}})
\ldots
\varphi(\Lambda_{B_{r}})
$$
where $\Lambda_{[s,t)}=\Lambda_{t}-\Lambda_{s}$ and 
$B_{1}, \ldots , B_{r}$ are the blocks of $R(\vec{k}, \vec{\sigma})$, 
with $R$ being 
the partition associated with the tuple $(l_{1}, \ldots , l_{n})$.
This fact follows from the continuous
tensor product decomposition of $\Gamma({\cal H})$ with respect
to time and the fact that 
all summands of $\Lambda_{I(l)}^{(k, \sigma)}$ have the form
$a(I(l))p(I(l),\sigma)$, where $p(I(l),\sigma)$ plays the role of 
${\bf P}(l,\sigma)$ in
(3.1), whereas $a(I(l))$ is an elementary tensor which has units
at all sites associated with $I(m)$'s for $m\neq l$. Note that they are not
filtered random variables in the sense of Definition 3.1, however (8.1)
still holds.
Moreover, 
the distribution $\phi_{I}=\phi|{\cal C}_{I}$, where ${\cal
C}_{I}$ is the unital *-algebra generated by
$\Lambda_{I}^{(k,\sigma)}$, $k\in {\bf N}$, $\sigma\in {\cal P}({\bf N})$,
depends only on 
the Lebesgue measure $\lambda (I)$ of $I$ 
since every expectation is in fact a polynomial
in the lenght of $I$. This can be seen by using Proposition 6.1.

In view of Lemma 8.2, it suffices to show that
$$
\lim_{N\rightarrow \infty}
\varphi(\Lambda_{t/N}^{(k,\sigma_{1})}\ldots \Lambda_{t/N}^{(k,\sigma_{n})})=t
$$
for any $k\in {\bf N}$ and $\sigma_{1}, \ldots , \sigma_{n}\in {\cal P}({\bf N})$.

Looking at the action of the fundamental filtered operators on the
finite particle domain (Proposition 6.1), we can see that 
in the considered expectation 
each creation-annihilation pair as well as each time operator produce 
$t$, whereas each number operator produces an integer. Therefore,
we obtain
$$
\varphi(\Lambda_{t}^{(k,\sigma)})=\varphi(A_{t}^{(0,\sigma)})=t,
$$
$$
\varphi(\Lambda_{t}^{(k,\sigma_{1})}\Lambda_{t}^{(k,\sigma_{2})})=
\varphi(A_{t}^{(k,\sigma_{1})}A_{t}^{(k,\sigma_{2})*})+o(t)
$$
$$
\varphi(\Lambda_{t}^{(k,\sigma_{1})}\ldots \Lambda_{t}^{(k,\sigma_{n})})
=
\varphi(A_{t}^{(k,\sigma_{1})}A_{t}^{(k,\sigma_{2})\circ}\ldots 
A_{t}^{(k,\sigma_{n-1})\circ}
A_{t}^{(k,\sigma_{n})*})	+
o(t)
$$
$$
=t+o(t)
$$
for $n>2$ and any $k\in {\bf N}$, $\sigma_{1}, \ldots , \sigma_{n}\in {\cal P}({\bf N})$. 
Hence 
$$
\lim_{N\rightarrow \infty}
N\varphi(\Lambda_{t/N}^{(k,\sigma_{1})}\ldots \Lambda_{t/N}^{(k,\sigma_{n})})=t
$$
which enables us to use Lemma 8.4 and obtain the desired form of
the expectation.\hfill $\Box$\\
\indent{\par}
The above theorem gives the combinatorics of the
filtered Poisson noise, defined by (8.1).
As it contains infinitely many colors and filters, this combinatorics
involves multivariate expectations (when speaking of a $1$-dimensional noise
we mean one ``type'' of operator, although it has infinitely many
``copies''). 
It can be noted that if $k\in \sigma$, then $\Lambda^{(k,\sigma)}_{t}$ 
for fixed $k$ and $\sigma$ gives classical Poisson white noise and 
if $k\notin \sigma$, then $\Lambda^{(k,\sigma)}_{t}$ gives Boolean (or, 1-free) 
Poisson white noise (cf. Theorem 5.4). 
Moreover, $m$-free and free
Poisson white noises are obtained from linear combinations of the same
type as in (3.3) and (5.4). In general,
also on the level of white
noise, filtered Gaussian and Poisson's white noises are also the building
blocks of other Gaussian and Poisson's white noises since the latter can be 
obtained from the former by addition or strong limits.
\\[10pt]
\myownsection
\begin{center}
{\sc 9. A free Fock space decomposition of $\Gamma({\cal H})$}
\end{center}
In this section we embed the free and $m$-free Fock spaces over ${\cal K}$, 
denoted by ${\cal F}({\cal K})$, ${\cal F}^{(m)}({\cal K})$, $m\in {\bf N}$,
respectively, in the multiple symmetric Fock space
$\Gamma ({\cal H})$, where ${\cal H}={\cal K}\otimes {\cal G}$, and
extend the $m$-free and free creation
and annihilation operators to bounded operators on $\Gamma({\cal H})$.
We assume that ${\cal K}=L^{2}({\bf R}_{+})$.

Let us introduce the following linear combinations of
filtered creation and annihilation operators, respectively:
\begin{eqnarray}
\label{9.1}
l^{(m)*}(f)&=&\sum_{k=1}^{m}(a^{(k)*}(f \otimes e_{k})-
a^{(k-1)*}(f\otimes e_{k}))\\
\label{9.2}
l^{(m)}(f)&=&\sum_{k=1}^{m}(a^{(k)}(f \otimes e_{k})-
a^{(k-1)}(f\otimes e_{k}))
\end{eqnarray}
where $m\in {\bf N}$. 
We will call $l^{(m)*}(f)$, $l^{(m)}(f)$, the {\it extended
$m$-free creation and annihilation operators}, respectively.
In order to compare them with the $m$-free creation
and annihilation operators $a^{(m)*}(f)$, $a^{*(m)}(f)$
introduced in [4], let us recall the definition of the latter.

First, the $m$-{\it free Fock space}
over ${\cal K}$ is the truncation of order $m$ of the free Fock
space, namely
$$
{\cal F}^{(m)}({\cal K})=
{\bf C}\omega_{m}\oplus \bigoplus_{k=1}^m {\cal K}^{\otimes k}
$$
where $\omega_{m}$ is the vacuum unit vector,
with the canonical scalar product.
The $m$-{\it free creation operators} are then given by
$$
a^{(m)*}(f):\; {\cal F}^{(m)}({\cal K})\rightarrow {\cal F}^{(m)}({\cal K})
$$
$$
a^{(m)*}(f)\; f_{1}\otimes \ldots \otimes f_{n} =
\left\{
\begin{array}{cll}
f\otimes f_{1}\otimes \ldots \otimes f_{n}& {\rm if} & 1\leq n < m\\
0& {\rm if} & n=m
\end{array}
\right.
$$
with $a^{(m)*}(f)\omega_{m}=f$ and the $m$-{\it free
annihilation operators}
$$
a^{(m)}(f):\; {\cal F}^{(m)}({\cal K})\rightarrow {\cal
F}^{(m)}({\cal K})
$$
$$
a^{(m)}(f)\; f_{1} \otimes \ldots \otimes f_{n} =
\langle f , f_{1} \rangle \; f_{2} \otimes \ldots \otimes f_{n}
$$
if $ 1 \leq n \leq m$ and $a^{(m)}(f)\omega_{m}=0$.
Note that $a^{(m)*}(f), a^{(m)}(f)$ are bounded on 
${\cal F}^{(m)}({\cal K})$ since they are
trunations of order $m$ of free creation and annihilation
operators $a^{*}(f)$, $a(f)$ on the free Fock space ${\cal F}({\cal K})$, 
respectively.
We will see below that $l^{(m)*}(f)$, $l^{(m)}(f)$ are bounded
extensions of $a^{(m)*}(f)$, $a^{(m)}(f)$, $f\in {\cal K}$, 
respectively, to $\Gamma({\cal H})$.

If we set $m=\infty$ in the formulas for extended $m$-free
creation and annihilation operators, we obtain operators 
which we denote $l^{*}(f)$ and $l(f)$, respectively
which will be called {\it extended free creation and annihilation
operators}. They, too, are bounded extensions of free creation and 
annihilation operators $a^{*}(f)$, $a(f)$, $f\in {\cal K}$,
to all of $\Gamma({\cal H})$, respectively. 

Thus, we identify two notations: $l^{(\infty)*}(f)\equiv l^{*}(f)$,
$l^{(\infty)}(f)\equiv l(f)$, $f\in {\cal K}$. In that
context we will understand that $P^{(\infty)}\equiv I$.
In general, in this section we will often assume for convenience
that $m\in {\bf N}^{*}={\bf N}\cup \{\infty\}$. 
However, certain results will be stated for $m=\infty$
separately in order to single out the free case. \\
\indent{\par}
{\it Remark.}
On the finite particle domain $\Gamma_{0}({\cal H})$,
spanned by $\Omega$ and vectors of the form
$$
(f_{1}\otimes e_{k_{1}})\circ \ldots \circ (f_{n}\otimes e_{k_{n}})
$$
where $f_{1}, \ldots , f_{n} \in {\cal K}$, 
$k_{1}\leq k_{2}\leq \ldots \leq k_{n}$, $n\in {\bf N}$,
the series given by (9.1)-(9.2) for $m=\infty$ are strongly convergent
since only a finite number of terms do not vanish
when acting on vectors of finite ``color support'' and thus give
well-defined operators with domains dense in $\Gamma({\cal H})$.
A similar feature 
was exhibited by the series representation of free random
variables obtained from the construction of the hierarchy of
freeness ([11], [5]).
We will see below that they have bounded extensions
to $\Gamma({\cal H})$.

Let us first determine the action of $m$-free creation and annihilation
operators on $\Gamma_{0}({\cal H})$. \\
\indent{\par}
{\sc Proposition 9.1.}
{\it Let $f,f_{1}, \ldots , f_{n}\in {\cal K}$ and 
$k_{1}\leq k_{2}\leq \ldots
\leq k_{n}$, $m\in {\bf N}^{*}$. Then}
$$
l^{(m)*}(f) \Omega = f\otimes e_{1} ,
$$
$$
l^{(m)*}(f) 
(f_{1}\otimes e_{k_{1}})\circ \ldots \circ (f_{n}\otimes e_{k_{n}})
$$
$$
=
\1_{\{m\geq k_{n}+1\}}
\sqrt{(n+1)}
(f_{1}\otimes e_{k_{1}})\circ \ldots \circ (f_{n}\otimes e_{k_{n}})
\circ (f\otimes e_{k_{n}+1}),
$$
$$
l^{(m)}(f)\Omega =0 ,
$$
$$
l^{(m)}(f)(f_{1}\otimes e_{k_{1}})\circ \ldots \circ (f_{n}\otimes e_{k_{n}})
$$
$$
=
\1_{\{m\geq k_{n}\}}
\frac{1}{\sqrt{n}}
\langle f , f_{n} \rangle
\delta_{k_{n},k_{n-1}+1} 
(f_{1}\otimes e_{k_{1}})\circ \ldots 
\circ (f_{n-1}\otimes e_{k_{n-1}})
$$
{\it where it is understood that $k_{0}=0$.}\\[5pt]
{\it Proof.}
Note that 
$$
a^{(k)*}(f \otimes e_{k})-a^{(k-1)*}(f\otimes e_{k})
=
a^{*}(f\otimes e_{k})P^{[k-1]}
$$
$$
a^{(k)}(f \otimes e_{k})-a^{(k-1)}(f\otimes e_{k})
=
P^{[k-1]}a(f\otimes e_{k})
$$
where 
$$
P^{[k-1]} 
(f_{1}\otimes e_{k_{1}})\circ \ldots \circ (f_{n}\otimes
e_{k_{n}}) 
=\delta_{k-1,k_{n}}
(f_{1}\otimes e_{k_{1}})\circ \ldots \circ (f_{n}\otimes
e_{k_{n}}) 
$$
since $P^{[k-1]}=P^{(k)}-P^{(k-1)}$ and $k_{j}\leq k_{n}$ for 
all $j=1, \ldots , n$. Thus,
$$
l^{(m)*}(f)
(f_{1}\otimes e_{k_{1}})\circ \ldots \circ (f_{n}\otimes e_{k_{n}})
$$
$$
=\sum_{k=1}^{m}a^{*}(f\otimes e_{k})P^{[k-1]}
(f_{1}\otimes e_{k_{1}})\circ \ldots \circ (f_{n}\otimes e_{k_{n}})
$$
$$
=\sum_{k=1}^{m}
\delta_{k-1,k_{n}}
a^{*}(f\otimes e_{k})
(f_{1}\otimes e_{k_{1}})\circ \ldots \circ (f_{n}\otimes e_{k_{n}})
$$
$$
=\1_{\{m\geq k_{n}+1\}}
a^{*}(f\otimes e_{k_{n}+1})
(f_{1}\otimes e_{k_{1}})\circ \ldots \circ (f_{n}\otimes e_{k_{n}})
$$
$$
=\1_{\{m\geq k_{n}+1\}}
\sqrt{n+1}
(f_{1}\otimes e_{k_{1}})\circ \ldots \circ (f_{n}\otimes e_{k_{n}})
\circ (f\otimes e_{k_{n}+1}).
$$
Next, if $n>1$, then
$$
l^{(m)}(f)
(f_{1}\otimes e_{k_{1}})\circ \ldots \circ (f_{n}\otimes e_{k_{n}})
$$
$$
=\sum_{k=1}^{m}P^{[k-1]}a(f\otimes e_{k})
(f_{1}\otimes e_{k_{1}})\circ \ldots \circ (f_{n}\otimes e_{k_{n}})
$$
$$
=\frac{1}{\sqrt{n}}\sum_{k=1}^{m}\sum_{j=1}^{n}
\langle f\otimes e_{k} , f_{j}\otimes e_{k_{j}}\rangle
P^{[k-1]}
(f_{1}\otimes e_{k_{1}})\circ \ldots \circ (f_{j}\breve{\otimes} e_{k_{j}})
\circ \ldots \circ (f_{n}\otimes e_{k_{n}})
$$
$$
=\frac{1}{\sqrt{n}}\sum_{k=1}^{m}\sum_{j=1}^{n-1}
\langle f, f_{j}\rangle
\delta_{k,k_{j}}
\delta_{k-1,k_{n}}
(f_{1}\otimes e_{k_{1}})\circ \ldots \circ (f_{j}\breve{\otimes} e_{k_{j}})
\circ \ldots \circ (f_{n}\otimes e_{k_{n}})
$$
$$
+
\frac{1}{\sqrt{n}}\sum_{k=1}^{m}
\langle f, f_{n}\rangle
\delta_{k,k_{n}}\delta_{k-1, k_{n-1}}
(f_{1}\otimes e_{k_{1}})\circ \ldots \circ (f_{n-1}\otimes e_{k_{n-1}})
$$
$$
=
\1_{\{m\geq k_{n}\}}
\frac{1}{\sqrt{n}}
\langle f, f_{n}\rangle
\delta_{k_{n}, k_{n-1}+1}
(f_{1}\otimes e_{k_{1}})\circ \ldots \circ (f_{n-1}\otimes e_{k_{n-1}})
$$
where the last equality follows again from the fact that $k_{j}\leq k_{n}$
for $j \leq n$, which makes the first sum vanish.
If $n=1$, then
$$
l^{(m)}(f)f_{1}\otimes e_{k_{1}}=
\sum_{k=1}^{m}\delta_{k, k_{1}}\langle f , f_{1} \rangle P^{[k-1]}\Omega
=\delta_{k_{1},1}\langle f , f_{1} \rangle  \Omega 
\equiv\delta_{k_{1},1}\langle f , f_{1} \rangle  
$$
The action on the vaccum vector is immediate.
\hfill $\Box$\\
\indent{\par}
{\sc Corollary 9.2.}
{\it In particular, if $m=\infty$, we obtain}
$$
l^{*}(f) \Omega = f\otimes e_{1} ,
$$
$$
l^{*}(f) 
(f_{1}\otimes e_{k_{1}})\circ \ldots \circ (f_{n}\otimes e_{k_{n}})
$$
$$
=
\sqrt{(n+1)}
(f_{1}\otimes e_{k_{1}})\circ \ldots \circ (f_{n}\otimes e_{k_{n}})
\circ (f\otimes e_{k_{n}+1}),
$$
$$
l(f)\Omega =0 ,
$$
$$
l(f)(f_{1}\otimes e_{k_{1}})\circ \ldots \circ (f_{n}\otimes e_{k_{n}})
$$
$$
= 
\frac{1}{\sqrt{n}}
\langle f , f_{n} \rangle
\delta_{k_{n},k_{n-1}+1} 
(f_{1}\otimes e_{k_{1}})\circ \ldots \circ (f_{n-1}\otimes e_{k_{n-1}}).
$$
\indent{\par}
{\sc Theorem 9.3.}
{\it For any $m\in {\bf N}^{*}$ and $f,g\in {\cal K}$, 
the operators $l^{(m)*}(f)$ and $l^{(m)}(f)$ have unique bounded
extensions to $\Gamma ({\cal H})$, 
are adjoints of each other, and satisfy the following relation:}
$$
l^{(m)}(g)l^{(m)*}(f)=\langle g ,f \rangle P^{(m)}.
$$
{\it Proof.}
By Proposition 9.1 we have
$$
l^{(m)}(f)l^{(m)*}(g)
(f_{1}\otimes e_{k_{1}})\circ \ldots \circ (f_{n}\otimes
e_{k_{n}}) 
$$
$$
=l^{(m)}(f)\1_{\{m\geq k_{n}+1\}}\sqrt{n+1}
(f_{1}\otimes e_{k_{1}})\circ \ldots \circ (f_{n}\otimes
e_{k_{n}}) 
\circ (g\otimes e_{k_{n}+1})
$$
$$
=\1_{\{m\geq k_{n}+1\}}
\langle f , g \rangle 
(f_{1}\otimes e_{k_{1}})\circ \ldots \circ (f_{n}\otimes
e_{k_{n}}) 
$$
where $k_{1}\leq k_{2} \leq \ldots \leq k_{n}$,
which proves that the relation holds on $\Gamma_{0}({\cal H})$. 
The proof of adjointness goes as follows.
$$
\langle
(f_{1}\otimes e_{k_{1}})\circ \ldots \circ (f_{n}\otimes
e_{k_{n}}) ,
l^{(m)*}(f)
(g_{1}\otimes e_{l_{1}})\circ \ldots \circ (g_{p}\otimes
e_{l_{p}})
\rangle 
$$
$$
=
\1_{\{m\geq l_{p}+1\}}
\sqrt{p+1}
(f_{1}\otimes e_{k_{1}})\circ \ldots \circ (f_{n}\otimes
e_{k_{n}}) ,
(g_{1}\otimes e_{l_{1}})\circ \ldots \circ (g_{p}\otimes
e_{l_{p}}) \circ (f\otimes e_{l_{p}+1})
\rangle
$$
$$
=
\1_{\{m\geq l_{p}+1\}}
\frac{\sqrt{n}}{n!}
\delta_{n,p+1}\delta_{k_{1},l_{1}}\ldots
\delta_{k_{n-1},l_{n-1}} \delta_{k_{n},l_{p}+1}
\langle f_{1}, g_{1} \rangle 
\ldots 
\langle f_{n-1}, g_{n-1} \rangle
\langle f_{n} , f \rangle .
$$
On the other hand,
$$
\langle
l^{(m)}(f)
(f_{1}\otimes e_{k_{1}})\circ \ldots \circ (f_{n}\otimes
e_{k_{n}}) ,
(g_{1}\otimes e_{l_{1}})\circ \ldots \circ (g_{p}\otimes
e_{l_{p}}) 
\rangle
$$
$$
=
\1_{\{m\geq k_{n}\}}
\frac{1}{\sqrt{n} (n-1)!}
\delta_{n-1,p}\delta_{k_{1},l_{1}}\ldots
\delta_{k_{n-1},l_{n-1}} \delta_{k_{n},k_{n-1}+1}
\langle f_{1} , g_{1} \rangle
\ldots
\langle f_{n-1} , g_{n-1} \rangle
\langle f_{n} , f \rangle
$$
where the following expression for the scalar product
$$
\langle
(f_{1}\otimes e_{k_{1}})\circ \ldots \circ (f_{n}\otimes
e_{k_{n}}) ,
(g_{1}\otimes e_{l_{1}})\circ \ldots \circ (g_{p}\otimes
e_{l_{p}}) 
\rangle
$$
$$
=
\frac{1}{n!}
\delta_{n,p}
\delta_{k_{1},l_{1}}
\ldots
\delta_{k_{n}, l_{n}} 
\langle f_{1}, g_{1} \rangle
\ldots 
\langle f_{n}, g_{n} \rangle .
$$
is obtained from the canonical scalar product on $\Gamma({\cal H})$.
Therefore, we have
$$
\langle l^{(m)}(f)x, y \rangle
=
\langle  x , l^{(m)*}(f) y \rangle
$$
for $x,y \in \Gamma_{0} ({\cal H})$. Now, note that
$$
\|l^{(m)*}(f)(f_{1}\otimes e_{k_{1}})
\circ \ldots \circ (f_{n}\otimes e_{k_{n}})\|^{2}=
\1_{\{m\geq k_{n}+1\}}
\|f\|^{2}
\|f_{1}\otimes e_{k_{1}})\circ \ldots \circ (f_{n}\otimes e_{k_{n}})\|^{2}
$$
hence $l^{(m)*}(f)$ has a unique bounded extension to  
$\Gamma({\cal H})$
of norm $\|l^{(m)*}(f)\|=\|f\|$ and thus
the annihilation operator $l^{(m)}(f)$ has also a unique bounded extension to 
$\Gamma({\cal H})$ of norm $\|l^{(m)}(f)\|=\|f\|$. \hfill $\Box$\\
\indent{\par}
Acting with the $m$-free creation and annihilation
operators on $\Omega$, $m\in {\bf N}$, and taking the closure,
we recover a subspace isomorphic 
to the $m$-free Fock space ${\cal F}^{(m)}({\cal K})$.
Thus, denote by $\widetilde{\cal F}^{(m)}({\cal K})$ the
closure of the space 
$\widetilde{\cal F}_{0}^{(m)}({\cal K})$ spanned by 
$\Omega$ and vectors of the form
$$
(f_{n}\otimes e_{1})\circ \ldots \circ (f_{1}\otimes e_{n})
$$ 
where $f_{1}, \ldots , f_{n}\in {\cal K}$, $1\leq n\leq m$ if $m$ is finite.
Similarly, denote by $\widetilde{\cal F}({\cal K})$ the closure
of $\widetilde{\cal F}_{0}({\cal K})$ spanned by vectors of the
above form with arbitrary $n\in {\bf N}$.
We obtain
$$
\langle (f_{n}\otimes e_{1}) \circ \ldots \circ (f_{1}\otimes
e_{n}) \rangle ,
(g_{m}\otimes e_{1}) \circ \ldots \circ 
(g_{1}\otimes e_{m}) \rangle
=\delta_{n,m}
\frac{1}{n!}
\langle f_{1}, g_{1} \rangle
\ldots
\langle f_{n} , g_{n} \rangle 
$$
by the orthogonality of $e_{1}, \ldots , e_{n}$.\\
\indent{\par}
{\sc Corollary 9.4.}
{\it The $m$-free Fock space ${\cal F}^{(m)}({\cal K})$ 
is isomorphic to $\widetilde{\cal F}^{(m)}({\cal K})$.
The free Fock space ${\cal F}({\cal K})$ 
is isomorphic to $\widetilde{\cal F}({\cal K})$.}\\[5pt]
{\it Proof.}
The unitary isomorphism from ${\cal F}_{0}({\cal K})$
to $\widetilde{\cal F}_{0}({\cal K})$ is given by 
$$
f_{1}\otimes \ldots \otimes f_{n}\rightarrow
\sqrt{n!}(f_{n}\otimes e_{1})\circ \ldots \circ (f_{1}\otimes e_{n})
$$
and thus extends uniquely to ${\cal F}({\cal K})$ (its restrictions
give the result for ${\cal F}^{(m)}({\cal K})$). \hfill $\Box$\\
\indent{\par}
Thus, for each $m\in {\bf N}$, we obtain the filtration
$$
\widetilde{\cal F}^{(1)}({\cal K})< \ldots <
\widetilde{\cal F}^{(m)}({\cal K}) <\ldots 
<\widetilde{\cal F}({\cal K})
$$
in which  
$\widetilde{\cal F}^{(m)}({\cal K})$ is an invariant subspace
for the $C^{*}$-algebra 
$$
{\cal C}^{(m)}=C^{*}\langle {\bf 1} , l^{(m)*}(f)|f\in {\cal K}\rangle
$$ 
and $\widetilde{\cal F}({\cal K})$ is an invariant subspace
for the $C^{*}$-algebra 
$$
{\cal C}=C^{*}\langle {\bf 1},l^{*}(f)|, f\in {\cal K} \rangle .
$$
Moreover, each $\widetilde{\cal F}({\cal K})$
is only one copy of the free Fock space in $\Gamma({\cal K})$
and it turns out that one can decompose $\Gamma({\cal H})$ into
a countable direct sum of subspaces isomorphic to
to the free Fock space and invariant under ${\cal C}$.
In the sequel we will concentrate on this decomposition,
in other words on what is ``between'' $\widetilde{\cal F}({\cal K})$
and $\Gamma({\cal H})$.

In order to determine this, we need to take a closer look 
at the kernel of the annihilation operators. 
Let $\{d_{n}\}_{n=1}^{\infty}$ be an orthonormal basis
in ${\cal K}$. 
Note that the set  
consisting of $\Omega$ and vectors of the form
$$
(d_{i_{1}}\otimes e_{k_{1}})\circ \ldots \circ (d_{i_{n}}\otimes e_{k_{n}})
$$
where $k_{1}\leq k_{2}\leq \ldots \leq k_{n}$ and $i_{r}\leq i_{r+1}$
whenever $k_{i_{r}}=k_{i_{r+1}}$, 
is an orthogonal basis in $\Gamma({\cal H})$ (the ordering
of indices is used for convenience, which is possible due to the
fact that the product is symmetrized).
Denote by $\widehat{\cal D}$ the subset of this basis 
consisting of $\Omega$ and vectors of the above form for which
$k_{n}\neq k_{n-1}+1$, i.e. 
the last two vectors are of identical colors or their
colors differ by more than $1$ if $n>1$,
and the last color is not equal to $1$ if $n=1$.
By normalizing the vectors from $\widehat{\cal D}$ we get 
$$
{\cal D}=\{x/\|x\|\; | x\in \widehat{\cal D}\}
$$
which is an orthonormal set. Let 
${\cal D}^{(m)}= {\cal D}\cap \Gamma^{m+1)}$, where
$\Gamma^{m+1)}= \Gamma({\cal H}^{(m+1)})$. We understand that
${\cal D}^{(\infty)}={\cal D}$. \\
\indent{\par}
{\sc Proposition 9.5.}
{\it ${\cal D}^{(m)}\subseteq {\rm ker}\;l^{(m)}(f)$
for any $m\in {\bf N}^{*}$ and $f\in {\cal K}$.}\\[5pt]
{\it Proof.}
This follows from Proposition  9.1 due to the presence of
$\delta_{k_{n}, k_{n-1}+1}$ on the right-hand side of the formula
for the annihilation operators.\hfill $\Box$\\
\indent{\par}
{\sc Proposition 9.6.}
{\it Let $m\in {\bf N}^{*}$. Then}
$$
\sum_{s=1}^{\infty}l^{(m)*}(d_{s})l^{(m)}(d_{s})=
I - P_{[{\cal D}^{(m)}]\oplus \Gamma^{(m}}
$$
{\it where $P_{[D^{(m)}]\oplus \Gamma^{(m}}$ is the projection onto
$[D^{(m)}]\oplus \Gamma^{(m}$ and 
$\Gamma^{(m}=\Gamma({\cal H}\ominus{\cal H}^{(m+1)})$. In particular,}
$$
\sum_{s=1}^{\infty}l^{*}(d_{s})l(d_{s})=
I - P_{[{\cal D}]}
$$
{\it thus ${\cal C}\cong {\cal O}_{\infty}$, where ${\cal O}_{\infty}$
is the Cuntz algebra.}\\[5pt]
{\it Proof.}
It can be seen from Theorem 9.3 that 
$$
l^{(m)*}(d_{s}) l^{(m)}(d_{s})=Q_{s},
$$
where $Q_{s}$ is the projection onto the subspace
spanned by vectors of the form
$$
(d_{s_{1}}\otimes e_{k_{1}})
\circ \ldots \circ
(d_{s_{n-1}}\otimes e_{k_{n-1}})
\circ
(d_{s}\otimes e_{k_{n}})
$$
where $k_{1}\leq \ldots \leq k_{n-1}=k_{n}-1<k_{n}\leq m$
(cf. [25]).
These subspaces are pairwise orthogonal and span the orthogonal 
complement of $[{\cal D}^{(m)}]\oplus \Gamma^{(m}$, which
proves the first formula. The second formula is just a special case
when $m=\infty$ and, together with Theorem 9.3, it implies that
the $C^{*}$-algebra generated by $l^{*}(f)$, $f\in {\cal K}$,
is isomorphic to the Cuntz algebra ${\cal O}_{\infty}$ 
since ${\cal K}$ is countably separable.
\hfill $\Box$\\
\indent{\par}
Let us introduce the following notation on $\Gamma_{0}({\cal H})$: 
$$
u\odot w=u_{1}\circ \ldots \circ u_{r} \circ w_{1} \circ \ldots \circ w_{n}
$$
where $u=u_{1}\circ \ldots \circ u_{r}$, 
$w= w_{1} \circ \ldots \circ w_{n}$. \\
\indent{\par}
{\sc Proposition 9.7.}
{\it If $x=x_{1}\circ \ldots \circ x_{r}$, $z=z_{1}\circ \ldots
\circ z_{r}$, $u=u_{1}\circ \ldots \circ u_{n}$,
$v=v_{1}\circ \ldots \circ v_{n}$, where 
$x_{i}, z_{i}\in {\cal H}_{1}$, $1\leq i \leq r$ and
$z_{j}, v_{j}\in {\cal H}_{2}$, $1\leq j \leq n$, and 
${\cal H}_{1},{\cal H}_{2}$ are two orthogonal 
subspaces of ${\cal H}$, then}
$$
\langle x\odot u , z\odot v \rangle
=
\frac{r!n!}{(r+n)!}
\langle x , z \rangle
\langle u , v \rangle .
$$
{\it Proof.}
Using the orthogonality of ${\cal H}_{1}$ and ${\cal H}_{2}$ and
the formula for the scalar product in $\Gamma({\cal H})$, we obtain
$$
\langle x\odot u , z\odot v \rangle
$$
$$
=\frac{1}{(r+n)!}\sum_{\sigma \in S_{r}}\sum_{\tau \in S_{n}}
\langle x_{1}, z_{\sigma(1)} \rangle
\ldots
\langle x_{r}, z_{\sigma(r)} \rangle
\langle u_{1}, v_{\tau(1)} \rangle
\ldots
\langle u_{n}, v_{\tau(n)} \rangle
$$
$$
=\frac{r!n!}{(r+n)!}
\langle x , z \rangle
\langle u , v \rangle .
$$
\hfill $\Box$\\ 
\indent{\par}
The $C^{*}$-algebra ${\cal C}$ 
is a $C^{*}$-subalgebra of ${\cal B}(\Gamma({\cal H}))$.
Denote the faithful representation of ${\cal C}$ on $\Gamma({\cal H})$
by $\pi$. 
Since $[{\cal C}x]$ is for each $x\in {\cal D}$
a closed subspace of $\Gamma({\cal H})$,
which is invariant under each operator $A$ in ${\cal C}$, the 
mapping $A\rightarrow A|[{\cal C}x]$ is a cyclic representation of 
${\cal C}$ on $[{\cal C}x]$ with cyclic vector $x$. 
Denote this representation by
$\pi_{x}$. We will show below that $\pi$ is a direct sum of cyclic
representations $\pi_{x}$, $x\in {\cal D}$.\\
\indent{\par}
{\sc Theorem 9.8.}
{\it The multiple symmetric Fock space has the direct sum decomposition}
$$
\Gamma({\cal H})=\bigoplus_{x \in {\cal D}}\;[{\cal C}x] 
$$
{\it where $[{\cal C}x]\cong {\cal F}({\cal K})$, according
to which}
$$
\pi=\bigoplus_{x\in {\cal D}} \pi_{x}
$$
{\it where $\pi_{x}\cong \rho$, and $\rho$ is the free Fock space
representation of ${\cal C}$.}\\[5pt]
{\it Proof.}
If $x=\hat{x}/\|x\|\in {\cal D}$, where $\hat{x}$ is of the form
$$
\hat{x}=(d_{1}\otimes e_{k_{1}})\circ \ldots \circ (d_{i_{r}}\otimes
e_{k_{r}})
$$
with $k_{r}=l$, then $[{\cal C}x]$ is the closed subspace of
$\Gamma({\cal H})$ spanned by vectors of the form
$$
x\odot 
(f_{n}\otimes e_{l+1})\circ \ldots \circ (f_{1}\otimes e_{l+n}) 
$$
where $f_{1}, \ldots , f_{n}\in {\cal K}$. Clearly, $[{\cal
C}x]$ is invariant under ${\cal C}$. Let us show that 
for each $x\in {\cal D}$, $[{\cal C}x]\cong {\cal F}({\cal K})$.

For that purpose, define the linear mapping
$$
U_{x}: \; {\cal F}_{0}({\cal K}) \rightarrow [{\cal C}x]
$$
by
$$
U_{x}(\omega)= x
$$
$$
U_{x}(f_{1}\otimes \ldots \otimes f_{n})=
\sqrt{\frac{(r+n)!}{r!}}
x \odot 
(f_{n}\otimes e_{l+1})\circ \ldots \circ (f_{1}\otimes e_{l+n}) .
$$
This mapping is scalar-product preserving since
$$
\langle U_{x}f_{1}\otimes \ldots \otimes f_{n},
U_{x}g_{1}\otimes \ldots \otimes g_{m}
\rangle
$$
$$
=\delta_{n,m}\frac{(r+n)!}{r!}
\langle 
x \odot 
(f_{n}\otimes e_{l+1})\circ \ldots \circ (f_{1}\otimes e_{l+n}),
x \odot 
(g_{n}\otimes e_{l+1})\circ \ldots \circ (g_{1}\otimes e_{l+n})
\rangle
$$
$$
=n!
\langle
(f_{n}\otimes e_{l+1})\circ \ldots \circ (f_{1}\otimes e_{l+n}), 
(g_{n}\otimes e_{l+1})\circ \ldots \circ (g_{1}\otimes e_{l+n})
\rangle
$$
$$
=\langle f_{1}, g_{1} \rangle
\ldots
\langle f_{n}, g_{n} \rangle ,
$$
and therefore has a unique extension to 
${\cal F}({\cal K})$. It is not hard to see that 
$[{\cal C}x]\perp [{\cal C}x']$ for $x\neq x'$ and that
$\Gamma({\cal H})$ is a direct sum of $[{\cal C}x]$ for
all $x\in {\cal D}$. 

It remains to be shown that $U_{x}$ intertwines between $\pi_{x}$
and the free Fock space 
representation $\rho$ of ${\cal C}$ on ${\cal F}({\cal K})$.
We have 
$$
\pi_{x}(l^{*}(f))U_{x}(\omega)=
\sqrt{r+1} x\odot (f\otimes e_{r+1})=U_{x}(f)=U_{x}\rho(l^{*}(f))\omega
$$
and 
$$
\pi_{x}(l^{*}(f))U_{x}(f_{1}\otimes \ldots \otimes f_{n})
$$
$$
=\sqrt{\frac{(r+n+1)!}{r!}}
x\odot 
(f_{n}\otimes e_{l+1})
\circ \ldots \circ
(f_{1}\otimes e_{l+n})
(f\otimes e_{l+n+1})
$$
$$
=U_{x}(f\otimes f_{1} \otimes \ldots \otimes f_{n})
$$
$$
=U_{x}\rho(l^{*}(f))(f_{1}\otimes \ldots \otimes f_{n})
$$
for any $f_{1}, \ldots , f_{n}, f\in {\cal K}$, $n\geq 1$.
Similarly, $\pi_{x}(l(f))U_{x}\Omega=l(f)x=0=U_{x}\rho(l(f))\Omega$
and
$$
\pi_{x}(l(f))U_{x}(f_{1}\otimes \ldots \otimes f_{n})
$$
$$
=\pi_{x}(l(f))
\sqrt{\frac{(r+n)!}{r!}}
x\odot (f_{n}\otimes e_{l+1})\circ \ldots \circ (f_{1}\otimes e_{l+n})
$$
$$
=\sqrt{\frac{(r+n-1)!}{r!}}
\langle f , f_{1} \rangle
x \odot 
(f_{n}\otimes e_{l+1})\circ \ldots \circ (f_{2}\otimes e_{l+n-1})
$$
$$
=\langle f , f_{1} \rangle 
U_{x}(f_{2}\otimes \ldots \otimes f_{n})
$$
$$
=U_{x}\rho(l(f))(f_{1}\otimes \ldots \otimes f_{n}).
$$
Therefore $\pi_{x}(a)U_{x}=U_{x}\rho(a)$ also for any $a\in {\cal C}$.
This finishes the proof.
\hfill $\Box$ \\
\indent{\par}
Let us finally define {\it extended $m$-free number 
operators}. Guided by the definitions of extended creation and 
annihilation operators, we set
$$
l^{(m)\circ}=\sum_{k=1}^{m}(a^{(k,k)\circ}-a^{(k,k-1)\circ})
$$
where $m\in {\bf N}^{*}$ and $a^{(k,r)\circ}=a^{(k,\sigma)\circ}$ for
$\sigma=\{1, \ldots , r-1\}$ and $a^{(k,\sigma)\circ}$ is given by (6.3).
Let us determine the action of extended number operators on the finite
particle domain.\\
\indent{\par}
{\sc Proposition 9.9.}
{\it Let $f_{1}, \ldots , f_{n}\in {\cal K}$, $k_{1}\leq \ldots \leq k_{n}$,
$m\in {\bf N}^{*}$. The finite particle domain $\Gamma_{0}({\cal H})$ is 
contained in the domains of extended $m$-free number operators and}
$$
l^{(m)\circ}
(f_{1}\otimes e_{k_{1}})\circ \ldots \circ (f_{n}\otimes e_{k_{n}})
$$
$$
=
\left\{
\begin{array}{cc}
N_{k_{n}}
(f_{1}\otimes e_{k_{1}})\circ \ldots \circ (f_{n}\otimes e_{k_{n}})&
{\rm if}\;\;k_{j}+1 =k_{n}\leq m\;\; {\rm for}\;\; j<n\\
0& {\rm otherwise}
\end{array}
\right.
$$
{\it where $N_{k}=\#\{i|k_{i}=k\}$.}\\[5pt]
{\it Proof.}
We have
$$
l^{(m)\circ}
(f_{1}\otimes e_{k_{1}})\circ \ldots \circ (f_{n}\otimes e_{k_{n}})
$$
$$
=\sum_{k=1}^{m}a^{(k)\circ}(P^{(\{1, \ldots ,k\})}-P^{(\{1, \ldots , k-2,k\})}
(f_{1}\otimes e_{k_{1}})\circ \ldots \circ (f_{n}\otimes e_{k_{n}})
$$
$$
=
\left\{
\begin{array}{cc}
a^{(k_{n})\circ}
(f_{1}\otimes e_{k_{1}})\circ \ldots \circ (f_{n}\otimes e_{k_{n}})&
{\rm if}\;\; k_{j}+1=k_{n}\leq m\;\;{\rm for}\;\; j<n\\
0& {\rm otherwise}
\end{array}
\right.
$$
$$
=
\left\{
\begin{array}{cc}
N_{k_{n}}
(f_{1}\otimes e_{k_{1}})\circ \ldots \circ (f_{n}\otimes e_{k_{n}})&
{\rm if}\;\; k_{j}+1=k_{n}\leq m\;\;{\rm for}\;\; j<n\\
0& {\rm otherwise}
\end{array}
\right.
$$
This ends the proof. \hfill $\Box$\\
\indent{\par}
In other words, $l^{(m)\circ}$ ``counts'' particles of the 
highest color $k_{n}$ if that one is smaller or equal to $m$ 
and the second highest color is equal to $k_{n}-1$. Otherwise,
the extended free number operator gives zero. In particular, 
on $\widetilde{\cal F}_{0}^{(m)}({\cal K})$ we obtain
$$
l^{(m)\circ}\Omega =0 
$$
$$
l^{(m)\circ}(f_{1}\otimes e_{n})\circ \ldots \circ (f_{n}\otimes e_{1})
=(f_{1}\otimes e_{n})\circ \ldots \circ (f_{n}\otimes e_{1})
$$
for $1 \leq n\leq m$.

It can be seen that, contrary
to the case of extended $m$-free creation and annihilation operators,
the operators $l^{(m)\circ}$ are not bounded on $\Gamma({\cal H})$.
Clearly, they are bounded on $\widetilde{\cal F}^{(m)}({\cal H})$
and, in fact, it can be shown that they are
bounded on $[{\cal C}x]$ for each $x\in {\cal D}$.\\[10pt]
\begin{center}
{\sc Acknowledgements}\\[10pt]
\end{center}
I would like to thank Piotr Hajac for stimulating discussions.\\[10pt]
\begin{center}
{\sc References}\\[10pt]
\end{center}
[1] {\sc L.~Accardi, M.~Sch\"{u}rmann, W.~von Waldenfels},
Quantum independent increment processes on superalgebra,
{\it Math. Z.} {\bf 198}, 451-477 (1988).\\[3pt] 
[2] {\sc A.~M.~Cockroft, R.~L.~Hudson}, Quantum mechanical Wiener 
processes, {\it J. Mult. Anal.} {\bf 7} (1977), 107-124.\\[3pt]
[3] {\sc A.~B.~Ghorbal, M.~Sch\"{u}rmann}, On the algebraic formulation 
of non-commu\-ta\-tive probability theory, preprint, Universite Henri
Poincare, 1999.\\[3pt]
[4] {\sc U.~Franz, R.~Lenczewski},
Limit theorems for the hierarchy of freeness, 
{\it Prob. Math. Stat.} {\bf 19} (1999), 23-41.\\[3pt]
[5] {\sc U.~Franz, R.~Lenczewski, M.~Sch\"{u}rmann}, 
The GNS construction for the hierarchy of freeness, Preprint No. 9/98,
Wroclaw University of Technology.\\[3pt]
[6] {\sc R.L.~Hudson, K.R.~Parthasarathy},
Quantum Ito's formula and stochastic evolution, 
{\it Commun.~Math.~Phys.} {\bf 93}, 301-323 (1984).\\[3pt]
[7] {\sc R.L.~Hudson, K.~R.~Parthasarathy}, Unification of
Fermion and Boson stochastic calculus, {\it Commun.~Math.~Phys.}
{\bf 104} (1986), 457-470.\\[3pt]
[8] {\sc B.~K\"{u}mmerer, R.~Speicher},
Stochastic integration on the Cuntz Algebra ${\cal O}_{\infty}$,
{\it J.~Funct.~Anal.} {\bf 103} (1992), 372-408.\\[3pt]
[9] {\sc R.~Lenczewski}, Quantum random walk for $U_{q}(su(2))$ 
and a new example of quantum noise, {\it J.~Math.~Phys} {\bf 37}
(1996), 2260-2278.\\[3pt]
[10] {\sc R.~Lenczewski}, Addition of independent variables in quantum 
groups, {\it Rev. Math. Phys.} {\bf 6} (1994), 135-147.\\[3pt]
[11] {\sc R.~Lenczewski}, Unification of independence in
quantum probability, {\it Inf.~Dim. Anal.~Quant.~Probab. \& Rel.~Top.}
{\bf 1} (1998), 383-405.\\[3pt]
[12] {\sc R.~Lenczewski}, Filtered stochastic calculus,
{\it Inf.~Dim.~Anal. Quant. Probab. \& Rel.~Top.}, to appear.\\[3pt]
[13] {\sc R.~Lenczewski, K.~Podgorski}, A $q$-analog of the 
quantum central limit theorem, {\it J.~Math.~Phys.} {\bf 33}
(1992), 2768-2778.\\[3pt]
[14] {\sc S.~Majid}, {\it Foundations of quantum groups}, Cambridge
University Press, 1995.\\[3pt]
[15] {\sc A.~Mohari, K.~B.~Sinha}, Quantum stochastic flows 
with infinite degrees of freedom and countable state Markov processes,
{\it Sankhya}, Ser.A, Part I (1990), 43-57.\\[3pt]
[16] {\sc K.~R.~Parthasarathy}, ``An Introduction to Quantum
Stochastic Calculus'', Birkha\"{a}user, Basel, 1992. \\[3pt]
[17] {\sc K.~R.~Parthasarathy, K.~B.~Sinha}, Unification of
quantum noise processes in Fock spaces,
Proc. Trento conference on Quantum Probability and Applications
(1989).\\[3pt]
[18] {\sc M.~Sch\"{u}rmann}, {\it White Noise on Bialgebras},
Lecture Notes in Math., Springer-Verlag, Berlin, 1993.\\[3pt]
[19] {\sc M.~Sch\"{u}rmann}, Non-commutative probability on algebraic
structures, Vol. XI (Oberwolfach, 1994), World Scientific, 1995, 
332-356.\\[3pt]
[20] {\sc M.~Sch\"{u}rmann}, Direct sums of tensor products and
non-commutative independence, {\it J. Funct. Anal.} {\bf 133}
(1995), 1-9.\\[3pt]
[21] {\sc R.~Speicher}, A new example of ``independence'' and ``white noise'',
{\it Probab.~Th. ~Rel.~Fields} {\bf 84} (1990), 141-159.\\[3pt]
[22] {\sc R.~Speicher, W.~von Waldenfels}, 
A general central limit theorem and
invariance principle, {\it in}: ``Quantum Probability and Related Topics'',
Vol.~IX, World Scientific, 1994, 371-387.\\[3pt]
[23] {\sc R.~Speicher, R.~Woroudi}, Boolean convolution, {\it
Fields Institute Commun.} {\bf 12} (1997), 267-279.\\[3pt]
[24] {\sc D.~Voiculescu}, Symmetries of some reduced free product 
${\cal C}^{*}$-algebras, {\it in} ``Operator Algebras and their 
Connections with Topology and Ergodic Theory'', Lecture
Notes in Math. 1132, Springer, Berlin, 1985, 556-588.\\[3pt]
[25] {\sc D.~Voiculescu, K.~Dykema, A.~Nica}, {\it Free Random
Variables}, CRM Monograph Series, AMS, Providence, 1992. 
\end{document}